\documentclass[12pt]{article}
\usepackage[english]{babel}
\usepackage{amsfonts,amsmath,amsxtra,amsthm,latexsym}
\usepackage{amssymb} 

\newtheorem{thm}{Theorem}[section]
 \newtheorem{cor}[thm]{Corollary}
 \newtheorem{lem}[thm]{Lemma}
 \newtheorem{prop}[thm]{Proposition}
 \theoremstyle{definition}
 \newtheorem{defn}{Definition}[section]
 \theoremstyle{remark}
 \newtheorem{rem}{Remark}[section]
 \numberwithin{equation}{section}

\DeclareMathOperator{\im}{Im}
\DeclareMathOperator{\re}{Re}

\DeclareMathOperator{\qarg}{qarg}

\def\R{\mathbb R}
\def\C{\mathbb C}
\def\N{\mathbb N}
\def\Z{\mathbb Z}
\def\D{\mathbb D}

\def\T{\mathbb T}

\def\B{\mathbb B}

\def\vphi{\varphi}
\def\la{\lambda}

\def\ep{\varepsilon}

\def\wt{\widetilde}
\def\dint{\displaystyle \int}

\def\i{\mathrm{i}}

\def\I{\mathcal{I}}
\def\Jf{\mathcal{J}}

\def\Ad{Ad}

\def\ChiCpl{\chi_{_{\scriptstyle \C_+}}}
\def\d{\mathrm{d}}
\def\k{\kappa}
\def\pa{\partial}

\def\r{\scriptstyle{r}}
\def\Sec{\mathrm{Sec}}
\def\Ext{\mathrm{Ext}}
\def\Exts{\mathrm{Ext}_\mathrm{step}}
\def\modn{\hspace{-10pt}\mod}

\def\wlim{\mathrm{w}^* \!\! - \! \lim}

\begin{document}
\title{Optimization of quasi-normal eigenvalues for 1-D wave equations in inhomogeneous media; description of optimal structures}
\author{}
\date{}
\maketitle
{\center {\large Illya M. Karabash } $^{\text{a,*}}$\\[4mm]

{\small $^{\text{a,*}}$  Institute of Applied Mathematics and Mechanics of NAS of Ukraine,

 R. Luxemburg str. 74, Donetsk, 83114, Ukraine. \ \ Tel.: +38 062 3110830. \ \ Fax: +38 062 3110285

E-mail addresses: i.m.karabash@gmail.com, karabashi@mail.ru

$^*$Corresponding author

}}

\begin{abstract}
The paper is devoted to optimization of resonances associated with
1-D wave equations in inhomogeneous media. The medium's structure is
represented by a nonnegative function $B$. The problem is to design
for a given $\alpha \in \R$ a medium that generates a resonance on
the line $\alpha + \i \R$ with a minimal possible modulus of the
imaginary part. We consider an admissible family of mediums that
arises in a problem of optimal design for photonic crystals. This
admissible family is defined by the constraints $0\leq b_1 \leq B
(x) \leq b_2$ with certain constants $b_{1,2}$. The paper gives an
accurate definition of optimal structures that ensures their
existence. We prove that optimal structures are piecewise constant
functions taking only two extreme possible values $b_1$ and $b_2$.
This result explains an effect recently observed in numerical
experiments. Then we show that intervals of constancy of an optimal
structure are tied to the phase of the corresponding resonant mode
and write this connection as a nonlinear eigenvalue problem.
\end{abstract}

\quad\\
MSC-classes: 49R05, 78M50, 35P25, 47N50, 47A55

\section{Introduction}

Recently the increasing interest in loss mechanisms of structured optical and mechanical systems has
 given rise to spectral optimization problems for dissipative models involving wave equations in inhomogeneous media,
 see e.g. \cite{AASN03,LSV03} and references therein.
The question is how to design an inhomogeneous medium with very low
radiative loss in a given frequency range. The radiative loss of
energy is closely connected to imaginary parts of eigenvalues of the
corresponding non-self-adjoint operator, see e.g. \cite{CZ95}. In
the paper these eigenvalues are called quasi-(normal) eigenvalues.
Naively, the closer quasi-eigenvalues to the real axis $\R$, the
less the radiative loss. In the recent numerical simulations
\cite{KS08,HBKW08} motivated by optimal design problems for photonic
crystals, the medium was modified by iterative methods with the
purpose to move a particular quasi-eigenvalue closer to $\R$.

An analytic background for spectral optimization problems involving
non-self-adjoint operators is not well developed. One of the
features that make non-self-adjoint spectral optimization problems
so different from self-adjoint ones is appearance of eigenvalues
with algebraic multiplicity greater than geometric multiplicity.
This leads to a much more complex perturbation theory for
eigenvalues.

The goal of the present paper is to study quasi-eigenvalue optimization problems analytically for
a 1-D model of a photonic
crystal with dissipation at one end.
To achieve this aim, we give a rigorous treatment of multiple eigenvalues and their perturbations.

More precisely, the paper is concerned with the eigenvalue problem
\begin{eqnarray}
y''(x) + \k^2 B(x) y (x) & = & 0 , \quad 0<x<1, \label{e ep}\\
y' (0) & =  & 0 , \label{e bc0}\\
y (1) - \i y' (1)/\k  & = & 0,
\label{e bc1}
\end{eqnarray}
where $\k$ is an eigen-parameter and the function $B$
satisfies
\begin{equation} \label{e a<b<bet}
 B\in L^{\infty} (0,1) , \qquad \ 0 \leq b_1 \leq B(x) \leq b_2 < \infty  \quad \text{for a.a.} \ x \in (0,1) .
\end{equation}

The spectral problem (\ref{e ep})-(\ref{e bc1}) comes from the
Fourier method applied to the scalar 1-D wave equation
\begin{eqnarray}
B(x) \pa_t^2 u (x,t) - \pa_x^2 u (x,t) = 0, \quad 0<x<1, \quad t>0, \label{e we}
\end{eqnarray}
equipped with the boundary conditions
\begin{eqnarray}
\pa_x u (0,t) =0 , \qquad \pa_x u (1,t) + \pa_t u (1,t) = 0 \label{e
bc we} .
\end{eqnarray}
The condition $\pa_x u (1,t) + \pa_t u (1,t) = 0$ leads to the
$\k$-dependent boundary condition (\ref{e bc1}) and corresponds to
radiative loss of energy to the surrounding medium through the
endpoint $x=1$.  Problems (\ref{e ep})-(\ref{e bc1}) and (\ref{e
we})-(\ref{e bc we}) arise in a number of applications.
 Let us mention spectral problems for a 1-D photonic crystal (see e.g. \cite[equation (2.26)]{JJWM08})
and for an inhomogeneous string damped at one end (see e.g. \cite{A75,KN79,KN89}).
In the first case the physical meaning of the function $B$ is the relative permittivity,
in the second case $B$ is the density of the string. Note also that $b_1$ is always positive for optical models,
but in the theory of a string, $B$ is allowed to be $0$ on a set of positive measure \cite{KK68_II,KN89}.
For studies concerned with the spectral problem (\ref{e ep})-(\ref{e bc1}) and the closely
related Regge problem we refer to \cite{A75,KN79,KN89,CZ95,F97,GP97,Sh07} and references therein.

Eigen-parameters $\k \in \C$ such that (\ref{e ep})-(\ref{e bc1})
has a nonzero solution will be called
\emph{quasi}-\emph{eigenvalues}. The corresponding eigenfunctions
are called (quasi-normal) \emph{modes}. Several other names for $\k$
are used, sometimes in slightly different settings: dissipation
frequencies \cite{KN79,KN89}, resonances \cite{KS08,HBKW08},
quasi-normal levels (in the Physics literature).

The set of quasi-eigenvalues is denoted by $K (B)$.
Quasi-eigenvalues $\k$ correspond to monochromatic solutions
$e^{i\k t} \vphi (x,\k)$ of the problem (\ref{e we})-(\ref{e bc we}).
The real part $\alpha=\re \k $ of the
quasi-eigenvalue is the \emph{frequency} of the monochromatic
solution, the imaginary part $\beta = \im \k$ is always positive and
characterizes the \emph{rate of decay}.

The following properties of quasi-eigenvalues are important for the
present paper: $K (B)$ is a subset of $\C_+$ symmetric with respect
to $\i \R$, quasi-eigenvalues are isolated, $\infty$ is their only
possible accumulation point, see e.g. \cite{KN79,KN89,CZ95}.

Let us explain the spectral optimization problem for (\ref{e
ep})-(\ref{e bc1}). We take the abstract point of view that the
problem (\ref{e we})-(\ref{e bc we}) is a mathematical model for a
certain device with a \emph{structure} $B (x)$. Assumption
(\ref{e a<b<bet}) defines the family $\Ad$ of admissible
structures.

We assume that the device is operated in a particular frequency range $[\alpha_1,\alpha_2]$, $-\infty < \alpha_1 \le \alpha_2 <+ \infty$, and denote by
$\Ad_{[\alpha_1,\alpha_2]} $ the set of all structures $B \in \Ad$
such that there exists at least one quasi-eigenvalue in this frequency range, i.e., such that
$K (B) \cap \{ z \in \C_+ \ : \ \re z \in [\alpha_1,\alpha_2] \} \neq \emptyset$.

\emph{The optimization problem under investigation is}
\begin{eqnarray}
 \text{to find} \ \  B_{\min} \in \Ad_{[\alpha_1,\alpha_2]} \ \ 
 \text{such that} \ \  \Jf_{[\alpha_1,\alpha_2]} (B_{\min}) \leq \Jf_{[\alpha_1,\alpha_2]} ( B ) \ \
\text{for all} \  B \in \Ad_{[\alpha_1,\alpha_2]} , \label{e OptJf int}
\end{eqnarray}
where
\begin{equation} \label{e def J}
\Jf_{[\alpha_1,\alpha_2]} (B) \ := \ \inf \{ \im \k \ : \ \re \k \in [\alpha_1,\alpha_2] \text{ and } \k \in K (B)
 \}  . 
\end{equation}

It seems that the systematic study of eigenvalue's maximization and
minimization problems associated with self-adjoint elliptic
operators was initiated by M.G.~Krein \cite{K51}. In the results of
\cite{K51} concerned with 1-D self-adjoint problems, extremizers are
extreme points of admissible families. Krein also proved the same
effect for one 2-D optimization problems and conjectured for another
\cite[Sec.4.4]{K51}. While there exists an extensive literature on
spectral optimization associated with self-adjoint elliptic
operators (see \cite{CMcL90,BOYa04} and references therein), there
are a very few analytically accurate papers on non-self-adjoint
spectral optimization problems similar to (\ref{e OptJf int}). A
possible explanation for this fact is that, for self-adjoint
problems, eigenvalues move on the real line and do not have root
eigenfunctions of higher order. This leads to a relatively simple
statement of the optimization problem and to a relatively simple
perturbation theory. Quasi-eigenvalues' behavior is much more
complex.

The problems of maximization of the  decay rate and of the spectral
abscissa for (\ref{e we})-(\ref{e bc we}) was considered in
\cite{CZ95}, where existence of the optimal design was proved for a
certain class of admissible structures from the space $W^{2,2}_\R
[0,1]$ and several estimates on quasi-eigenvalues were obtained.

In mathematical modeling for photonic crystals, the relative permittivity $B$ is usually discontinuous.
That is why the admissible family  (\ref{e a<b<bet}) is a reasonable choice.
This admissible family with $b_1 =1 $ was used recently in \cite{KS08}, where a gradient ascent iterative procedure
for optimization of the quality factor $Q(\k ) = \frac{|\re \k|}{2 |\im \k|} $
for individual quasi-eigenvalues $\k$ has been developed.
Numerical computations of \cite{KS08} were done for the 1-D and 2-D scalar wave equations.
Another numerical paper \cite{HBKW08} is concentrated on the 1-D case, but include
 an additional coefficient $\sigma (x)$ into
the equation $\pa_x \sigma \pa_x y  + \k^2 B y  =  0 $ and considers
various admissible families and discretization techniques.
Simulations of \cite{HBKW08} were performed for the case $B \equiv
1$. It was noticed that for a problem with constraints $1 = \sigma_1
\leq \sigma(x) \leq \sigma_2 = 3$ the optimization procedure stopped
on a structure $\sigma$ taking only the extreme possible values
$\sigma_{1,2}$. In 1-D simulations of \cite{KS08}, the coefficient
$B$ also tends to be a piecewise constant function taking values
$b_{1,2}$ (see \cite[Fig. 1-2 and page 423]{KS08}). In the author's
opinion, figures 3 and 4 of \cite{KS08} suggest the same effect in
the 2-D case.

The main results of the present paper are collected in Section
\ref{s pr inf} (except Theorem \ref{t sw point}, which requires more
preliminaries). We adjust the existence of minimizer  proof of the
self-adjoint case \cite{K51} to prove that the set of all possible
quasi-eigenvalues $K (\Ad) := \bigcup_{B \in \Ad} K (B) $ is closed.
This easily implies that the minimum of $\Jf_{[\alpha_1,\alpha_2]}$
is archived and is positive whenever the domain of definition of $
\Jf_{[\alpha_1,\alpha_2]} $ is nonempty (see Corollary \ref{c
SolInf}). Then the problem (\ref{e OptJf int}) can be reduced to the
study of the case when $\alpha_1 = \alpha_2 = \alpha$. We introduce
the function $\I (\alpha)$, $\alpha \in \R$, as the minimum of the
functional $\Jf_{[\alpha,\alpha]}$. It is natural to call the
complex points on its graph $\{(\alpha, \I (\alpha))\}_{\alpha \in
\R} = \{\alpha+\i \I (\alpha)\}_{\alpha \in \R}$ \emph{optimal
quasi-eigenvalues}, and to call structures $B$ corresponding to $\k
= \alpha+\i \I (\alpha) $ \emph{optimal structures} (for the
frequency $\alpha$). (The author does not know whether it is
possible that there are non-equivalent optimal structures
corresponding to a certain frequency $\alpha \neq 0$). \emph{Optimal
modes} for the frequency $\alpha$ are eigenfunctions of (\ref{e
ep})-(\ref{e bc1}) with corresponding optimal $\k$.

Theorems \ref{t B in Exts} and \ref{t J alpha0} state that
\emph{optimal structures are piecewise constant functions taking
only values $b_1$ and $b_2$}. For $\alpha =0$, we find $\I (0)$ and
the corresponding optimal structure explicitly. The effect behind
Theorem \ref{t nonlin} is that \emph{the intervals where an optimal
structure $B$ takes the values $b_1$ or $b_2$ are connected with the
$\arg$-function of the corresponding optimal mode.} This connection
can be written in a form of a nonlinear eigenvalue problem. That is,
if $\k_0$ is an optimal quasi-eigenvalue, then the equation
\begin{equation} \label{e nonlin eq int}
y'' + \k_0^2 y  \left[ b_1 + (b_2 - b_1) \ChiCpl (y^2 ) \right] =  0
\ \ \text{a.e. on } \ (0,1) ,
\end{equation}
has a non-trivial solution $y_0$ satisfying boundary conditions
(\ref{e bc0})-(\ref{e bc1}). Here $\ChiCpl (z) := 1$ for $z \in \C_+
$, and $\ChiCpl  (z) := 0$ for $z \in  \C \setminus \C_+$. This
solution $y_0$ is an optimal mode corresponding to $\k_0$. The
optimal structure associated with $\k_0$ and $y_0$ is $B (x) = b_1 +
(b_2 - b_1) \ChiCpl \left(y_0^2 (x) \right)$.

The results of Section \ref{s MultEig} on perturbations of
quasi-eigenvalues are preparative for the proofs of Theorems \ref{t
B in Exts} and \ref{t J alpha0}. Perturbations for quasi-eigenvalues
of the Schr\"{o}dinger operator $- \Delta +V $ were studied, e.g.,
in \cite{GK71,AAD01}, for an abstract approach and more references
see \cite{A98}. However the proofs of our results require more
delicate information. Namely, we need analyticity of a
quasi-eigenvalue as a functional of $B$ and we perform a detailed
(though non-complete) study of these functionals in vicinity of
their singular points. Roughly speaking, these singularities
correspond to multiple quasi-eigenvalues. The proofs of the main
results are based on Lemma \ref{l 2par per}, which is essentially
concerned with two parameter perturbations of a multiple
quasi-eigenvalue. Note that the gradient algorithm meets obvious
difficulties when it encounters a multiple eigenvalue, see the
discussion in \cite[p. 425]{KS08}. An accurate treatment of a
multiple eigenvalue requires understanding of its splitting picture.
Proposition \ref{p per k mult} and Lemma \ref{l 2par per} provide a
part of this collision and splitting picture.

Section \ref{s pr t extr p} contains the proofs of Theorems \ref{t B
in Exts} and \ref{t J alpha0}. The proofs are based on the
perturbation results and on the detailed study of a special solution
$\vphi$ of (\ref{e ep}) singled out by $\vphi (0) =1$, $\vphi ' (0)
= 0$. In Section \ref{s sw pt}, we prove Theorem \ref{t nonlin} and
study the interplay of optimal structures $B (x)$ and phases $\arg
\vphi (x)$ of associated optimal modes. For the non-degenerate case
when $b_1 >0$ and $\alpha \neq 0$, this interplay is written in an
especially transparent form as Theorem \ref{t sw point}: \emph{there
exists $\omega \in [-\pi,\pi)$ such that $B$ changes its value from
$b_1$ to $b_2$ exactly when $\vphi^2 $ crosses the ray $e^{\i
\omega} \R_+$ and from $b_2$ to $b_1$ exactly when $\vphi^2 $
crosses the ray $e^{\i \omega } \R_-$.}


\textbf{Notation}. $\C_\pm = \{ z \in \C : \pm \im z >0 \}$, \ $\R_\pm = \{ x \in \R: \pm x >0 \}$,
 $\D_\epsilon (\zeta) := \{z \in \C : | z - \zeta | < \epsilon \}$,
$\T = \{ z \in \C : | z | = 1 \}$.
For $\xi_1,\xi_2$ such that $0<\xi_2-\xi_1<\pi $, $\Sec
(\xi_1,\xi_2)$ defines the sector (without zero)
\[
\Sec (\xi_1,\xi_2) := \{ \zeta \in \C\setminus{0} \, : \, \arg \zeta = \xi \ (\modn 2\pi)
 \text{ for certain } \xi \in [\xi_1, \xi_2] \} .
\]

$\chi_E $ is an indicator function of the set $E $, i.e., $\chi_E
(x) = 1$ when $x \in E$, and $\chi_E (x) = 0$ when $x \not
\in E$. 

Open balls in a normed space $V$ are denoted by $\B_\epsilon (v_0) := \{v \in V \, : \, \| v - v_0 \|_V < \epsilon \}$.
For $\Omega \subset V$ (including the case $V=\C$), $v_0 \in V$, and $z \in \C$, let $z\Omega +
v_0 := \{ zv + v_0 \, : \, v \in \Omega \}$. For a function $f$
defined on $\Omega \subset V$, $f(\Omega)$ is the image of $\Omega$.

$L_{\C (\R)}^p (0,1) $ are the Lebesgue spaces of functions with values in $\C$
(resp., $\R$);
\[
W_\C^{k,p} (0,1) := \{ y \in L_\C^p (0,1) \ : \ \pa_x^j y \in L_\C^p (0,1), \ \ 1 \leq j \leq k \}
\]
are Sobolev spaces with standard norms. The space of continuous
complex-valued functions with the uniform norm is denoted by $C
[0,1]$.

$\pa_x y$, $\pa_z \vphi (x,z;B)$, etc. denote (ordinary or partial) derivatives with respect to (w.r.t.) $x$, $z$, etc.;
$[\pa_B \vphi (x,z;B)] (B_\Delta) = \lim_{\zeta \to 0} \frac{\vphi (x,z;B + \zeta B_\Delta) - \vphi (x,z;B)}{\zeta} $ is the directional derivative of the functional $\vphi (x,z;\cdot)$ along the direction $B_\Delta \in L_\C^\infty (0,1) \setminus \{0\}$ at the point $B \in L_\C^\infty (0,1)$.

We write $z_1^{[n]} \asymp z_2^{[n]}$ as $n \to \infty$ if the sequences
$z_1^{[n]} / z_2^{[n]} $ and $z_2^{[n]} / z_1^{[n]}$ are bounded for $n$ large enough.

\section{Optimal structures, the definition and main results}
\label{s pr inf}

Recall that the set of quasi-eigenvalues corresponding to a structure $B$ (in short, quasi-eigenvalues of $B$)
is denoted by $K(B)$. It occurs  that $K (B)$ is the set of zeroes of the entire function
\begin{equation*} \label{e F}
F(z) = F (z; B ) := \vphi (1,z) - \i \pa_x \vphi (1,z) \, / \, z , \ \ z \in \C,
\end{equation*}
where $\vphi (x,z) = \vphi (x, z; B )$ is the solution of the initial value problem
\begin{equation*} \label{e phi}
\pa_x^2 y(x,z) = - z^2 \ B (x) \ y (x,z) , \ \ \ y (0,z) = 1, \ \ \ \pa_x y (0,z) = 0 .
\end{equation*}
It is obvious that all modes $y$ corresponding to
$\k \in K (B)$ are equal to $\vphi $ up to a multiplication by a constant.
So \emph{the geometric multiplicity} of any quasi-eigenvalue equals 1.
In the following, \emph{the multiplicity} of a quasi-eigenvalue means its \emph{algebraic multiplicity}.

\begin{defn} \label{d mult}
\emph{The multiplicity} of a quasi-eigenvalue is its multiplicity as a zero of the entire function $F (\cdot) $.
A quasi-eigenvalue is called \emph{simple} if its multiplicity is $1$.
\end{defn}

This is classical M.V. Keldysh's definition of multiplicity for eigenvalue problems
with an eigen-parameter in boundary conditions, see e.g. \cite[Sec. 1.2.2-3]{N69} and \cite{KN79,KN89}.

Each quasi-eigenvalue has finite multiplicity.
The set of quasi-eigenvalues $K(B)$ is always symmetric w.r.t. the imaginary axis $\i \R$, moreover,
the multiplicities of symmetric quasi-eigenvalues are the same.
Note that $K (B)$ may be empty, this is the case for $B \equiv 0$ and $B \equiv 1$.
These and other basic facts can be found in \cite{KN79,KN89} (see also \cite{CZ95}).

We consider the quasi-eigenvalue problem (\ref{e ep})-(\ref{e bc1}) over the following family of structures
\begin{eqnarray*} \label{e Adinf}
\Ad  :=
\{ B \in L_\R^{\infty} (0,1) \ : \ b_1 \leq B(x) \leq
b_2 \ \ \text{a.e.}  \} , \ \ 0 \leq b_1 \leq b_2  < \infty, \ \ b_2 >0 .
\end{eqnarray*}



\begin{prop} \label{p SetInf}
The set $K (\Ad) := \bigcup_{B \in \Ad} K (B) $ is closed and is a subset of $\C_+$.
\end{prop}

The first statement is proved in the next subsection, see Lemma \ref{l compomega}.
The statement $K (\Ad) \subset \C_+ $ follows from the well-known fact that
$K (B) \subset \C_+$ for any $B \in \Ad$, see e.g. \cite{KN79,KN89}.

Now we pass to an immediate corollary, which shows that if $K (\Ad) $ has at least one
$\k$ with frequency in the range $[\alpha_1,\alpha_2]$ ($\alpha_1,\alpha_2 \in \R$, $\alpha_1 < \alpha_2$),
then a minimizer $B_{\min}$ for the optimization problem (\ref{e OptJf int}) exists.

Recall that the functional $\Jf_{[\alpha_1,\alpha_2]} $ is defined by (\ref{e def J})
on the family
\begin{eqnarray*} 
\Ad_{[\alpha_1,\alpha_2]}
:=
\{ B \in \Ad \ : \  \re \k \in [\alpha_1,\alpha_2] \text{ for certain } \k \in K (B) \} .
\end{eqnarray*}



\begin{cor} \label{c SolInf}
Suppose $\Ad_{[\alpha_1,\alpha_2]} \neq \emptyset$ and consider the problem (\ref{e OptJf int}). Then:
\begin{itemize}
\item[(i)] The functional $\Jf_{[\alpha_1,\alpha_2]} $ achieves its minimum
$\I_{[\alpha_1,\alpha_2]} := \Jf_{[\alpha_1,\alpha_2]} (B_{\min})$ over $\Ad_{[\alpha_1,\alpha_2]} $.
\item[(ii)] The minimum $\I_{[\alpha_1,\alpha_2]}$ is positive.
\end{itemize}
\end{cor}

A simple way to check the condition $ \Ad_{[\alpha_1,\alpha_2]} \neq \emptyset $ is to consider
'constant' structures $B \equiv b \in [b_1,b_2]$ (here and below $B
\equiv b$ means $B(x)=b$ for a.a. $x \in [0,1]$). For them the
quasi-eigenvalues are well known.

\begin{prop}[see e.g. \cite{CZ95}] \label{p const struct}
Let $B \equiv b$ be a constant function with $b \geq 0$. Let $\{ \k_n \} = K (B)$ be the set of corresponding quasi-eigenvalues (taking multiplicities into account).
Then:
\item[(i)] If $b=0$ or $b=1$, then $K (B) = \emptyset$.
\item[(ii)] If $b \not \in \{ 0, 1 \}$, then
$ \k_n = \i \frac{1}{2\sqrt{b}} \log \left| \frac{\sqrt{b}+1}{\sqrt{b}-1} \right| +
\frac{\pi}{\sqrt{b}} \left\{ \begin{array}{ll}
n , & \text{if } b>1\\
n+1/2, & \text{if } b<1
\end{array} \right. , \ \ n \in \Z . $
\end{prop}

So, excluding the extreme case  $b_1 = b_2 = 1$, we can always ensure
$\Ad_{[\alpha_1,\alpha_2]} \neq \emptyset$ taking the frequency range
$[\alpha_1,\alpha_2]$ wide enough.

Consider the case $\alpha_1=\alpha_2 = \alpha $ and introduce the function:
\begin{equation} \label{e def I inf}
 \I  (\alpha) \ := \ \left\{ \begin{array}{l}
 +\infty  \ \ \text{ if } \Ad_{[\alpha,\alpha]} = \emptyset \\
  \min \{ \im \kappa \ : \ \re \kappa = \alpha \text{ and } \k \in K (\Ad) \}  \ \text{ if } \Ad_{[\alpha,\alpha]} \neq \emptyset  \\
\end{array} \right.
\end{equation}
By Corollary \ref{c SolInf}, $\I  (\alpha) > 0 $ for all $ \alpha \in \R$
and the minimum in (\ref{e def I inf}) is achieved whenever $\Ad_{[\alpha,\alpha]} \neq \emptyset$.
Obviously, the minimal value $\I_{[\alpha_1,\alpha_2]}$ of the functional $\Jf_{[\alpha_1,\alpha_2]} (\cdot)$  is given by
$ \min_{\alpha \in [\alpha_1,\alpha_2]} \I  (\alpha) $. So problem  (\ref{e OptJf int}) can be reduced to the study of the function
$\I  $ and the properties of structures corresponding to
quasi-eigenvalues of the form $\k=\alpha + \i \I (\alpha)$.
 It is natural to call such $\k$ and $B$ \emph{optimal} (they are optimal at least for a particular frequency $\alpha$).

\begin{defn}
Let $\I (\alpha) < \infty$ for certain $\alpha \in \R$.
Then:
\item[(i)] $\k=\alpha + \i \I (\alpha)$ is called \emph{an optimal quasi-eigenvalue} for the frequency $\alpha$,
\item[(ii)] a structure $B \in \Ad $ is called \emph{optimal} for the frequency $\alpha$ if
$\alpha + \i \I (\alpha) \in K (B)$.
\end{defn}

One can check that the set $\Ext $ of extreme points of $\Ad$ is
\[
\Ext = \left\{ B \in \Ad \, : \, B(x) \in \{ b_1, b_2 \} \text{ for a.a. } x \in [0,1] \right\} .
\]
Denote
\begin{equation} \label{e E12}
E_j (B) := \{ x \in [0,1] \, : \, B(x) = b_j \}, \ \ j=1,2.
\end{equation}
Recall that a function $B$ on $[0,1]$ is called \emph{a piecewise constant function} if there exists a partition
$0 = x_0 < x_1 < x_2 < \dots < x_{n} < x_{n+1} = 1$
 such that $B$ is constant on each interval $(x_{j-1} , x_{j}) $.
By $\Exts $ we denote the family of  piecewise constant functions that belong to $\Ext$ (more precisely, the family of corresponding classes of equivalence), i.e.,
\begin{eqnarray*}
\Exts = \{ B \in \Ext \, & : & \, E_1 (B), \, E_2 (B) \text{ are unions of a finite number of intervals} \\ & & \text{after possible correction on sets of zero measure} \} .
\end{eqnarray*}

\begin{thm} \label{t B in Exts}
Assume that $B$ is an optimal structure for a frequency $\alpha \in \R$.
Then $B \in \Exts$, i.e., $B$ is a piecewise constant function taking only values $b_1$ and $b_2$
(after possible correction on a set of measure zero).
\end{thm}

The proof is given in Section \ref{s pr t extr p}, see Corollary \ref{c B in Exts}, Proposition \ref{p k nonmin alpha0}, and
the proof of Theorem \ref{t J alpha0}.

For $\alpha=0$, we find the optimal quasi-eigenvalue and the corresponding structure explicitly.

\begin{thm} \label{t J alpha0}
(i) If $b_2 \leq 1$, then $\I (0) = + \infty$.\\[1mm]
(ii) If $b_2 > 1$, then $\I (0) = \frac{1}{2 \sqrt{b_2}} \log \frac{\sqrt{b_2} +1}{\sqrt{b_2} -1}$, and the only structure in $\Ad$ having the quasi-eigenvalue at $\i \I (0)$ is $B \equiv b_2$.
\end{thm}

The proof is given in Section \ref{ss case a=0}. We would like to
note that statement (i) is equivalent to the fact that $b_2 \leq 1$
implies $K (\Ad) \cap \i \R = \emptyset$. Under the additional assumption $0<b_1<b_2 <1$, this fact was
obtained in \cite[Theorem 4.2 (i)]{CZ95} (our proof is completely different).

In the general case $\alpha \in \R$, the intervals where an optimal structure $B$ takes the values $b_1$ or $b_2$ are
connected with the $\arg$-function of the corresponding mode.
This connection can be written in a concise way as a nonlinear eigenvalue problem.


Put $\ChiCpl  (z) := 1$ when $\im z > 0$, and $\ChiCpl  (z) := 0$ when $\im z \leq 0$.
Consider the nonlinear equation
\begin{equation} \label{e nonlin eq}
y''(x) + \k^2 y (x) \left[ b_1 + (b_2 - b_1) \ChiCpl \left(y^2 (x) \right) \right] =  0 \ \
\text{a.e. on } \ (0,1) .
\end{equation}

\begin{thm} \label{t nonlin}
Let $0\leq b_1 < b_2$. Let $\k$ be an optimal quasi-eigenvalue for a
frequency $\alpha \in \R$.
Then there exists a nonzero solution $y \in
W_\C^{2,\infty} [0,1]$ of the nonlinear boundary value problem
(\ref{e nonlin eq}), (\ref{e bc0}), (\ref{e bc1}). Moreover, $B (x)
= b_1 + (b_2 - b_1) \ChiCpl \left(y^2 (x) \right)$ is an optimal structure
for the frequency $\alpha$.
\end{thm}

Theorem \ref{t nonlin} is proved in Section \ref{s sw pt}, where the
connection of the rotation of $\vphi^2 (x,\k;B)$ around $0$ with
intervals of constancy of the optimal structure $B$ is explained in
details. Note that the solution $y$ of  Theorem \ref{t nonlin} is
also a mode of the original linear problem (\ref{e ep}), (\ref{e
bc0}), (\ref{e bc1}). So $y (\cdot) = c \vphi (\cdot , \k;B)$ with
some constant $c$. In Section \ref{s sw pt}, we find appropriate
constants $c = e^{\i \theta}$. Generally, $\theta \neq 0 \ (\, \modn
2 \pi)$ and $\vphi$ may be \emph{not} a solution of (\ref{e nonlin
eq}), (\ref{e bc0}), (\ref{e bc1}). However, in some cases $y$ can
be taken equal to $\vphi$.

\begin{rem}
(1) Some estimates on $\I (\alpha)$ from above can be easily
obtained from Proposition \ref{p const struct}. For a certain range
of $\alpha$ estimates on $\I (\alpha)$ from below can be obtained
using a solution of the direct spectral problem for strings of the
Krein-Nudelman class \cite{Ka12_optKN}.

(2) Theorem \ref{t nonlin} does \emph{not} state that $B$ is a
unique optimal structure for the frequency $\alpha$. The author does
not know whether it is possible that there are non-equivalent
optimal structures corresponding to certain $\alpha \neq 0$. For
quasi-eigenvalue optimization problems in classes of Krein strings
with total mass and statical moment constraints optimal structures
are unique for their $\alpha$ \cite{Ka12_optKN}.
\end{rem}

\subsection{Proof of Proposition \ref{p SetInf}. \label{ss proof infty}}

\begin{lem}[Integral form of (\ref{e ep})-(\ref{e bc1})] \label{l int}
A number $\k \in \C$ belongs to $K (B)$ if and only if there exists $y (x) \in C [0,1]$ such that
\begin{eqnarray}
 y (x) & = & 1 - \k^2 \ \dint_0^x (x-s) \ B (s) \ y (s) \ \d s ,
\quad 0 \leq x \leq 1,
\label{e int ep} \\
y (1) & + & \i \k \int_0^1 B(s) y (s) \d s = 0
 \ . \label{e int bc1}
\end{eqnarray}
If such $y$ exists, then $y(x)= \vphi (x,\k;B)$.
\end{lem}

\begin{proof}
Equality (\ref{e int ep}) holds exactly when $y (\cdot) = \vphi (\cdot,\k;B)$. Using (\ref{e int
ep}), one can derive (\ref{e int bc1}) from (\ref{e bc1}).
And vise versa, equalities (\ref{e int ep})-(\ref{e int bc1}) imply
that $\k \neq 0$. Indeed, if $\k=0$, one has $y(1)=0$, which
contradicts (\ref{e int ep}). Finally, for $\k \neq 0$, we can
rewrite (\ref{e int bc1}) as (\ref{e bc1}) using (\ref{e ep}).
\end{proof}

\begin{lem} \label{l Gron}
The mapping $(z,B) \to \vphi (\cdot, z; B)$ is bounded from $\C \times L_\C^\infty$ to $W_\C^{2,\infty} [0,1]$.
\end{lem}

\begin{proof}
If the families $\{ z_\omega \}$ and $\{ B_\nu \} $
are bounded subsets of $ \C$ and $L_\C^\infty (0,1)$, resp.,
all functions $\vphi (\cdot, z_\omega; B_\nu )$ satisfy the differential inequality
$| y'' (x) | \leq C_1 |y(x)| $.
In turn, this implies (e.g., via the
Gronwall-Bellman inequality applied to $|y|$)
the statement of the lemma.
\end{proof}


\begin{lem}[cf. \cite{K51} for the self-adjoint case]\label{l compomega}
Assume that there exist sequences $\{ \k_n \}_1^\infty \subset \C$ and $\{ B_n \}_1^\infty \subset \Ad $ such that
$\k_n \in K(B_n) $ and $\k_n \to \wt \k \in \C$ .
Then there exists $\wt B \in \Ad $ such that
$\wt \k \in K (\wt B )$. 
\end{lem}

\begin{proof}
By the sequential Banach--Alaoglu theorem, there exist $\wt B \in
\Ad $ and a subsequence $\{ B_{n_j} \}$  such that $\wlim B_{n_j} =
\wt B$  (in weak* topology of $L^\infty $). Since $\k_{n_j} \to \wt
\k$, Lemma \ref{l Gron} yields that the sequence $\{ \vphi
(\cdot,\k_{n_j}; B_{n_j} ) \}$ is bounded in $W_\C^{2,\infty}
[0,1]$. The embedding $W_\C^{2,\infty} [0,1] \Subset C [0,1]$ is
compact, so there exists a subsequence $\{ m_j\}$ of $\{ n_j \}$
such that $ \{ \vphi (\cdot,\k_{m_j}; B_{m_j}  ) \} $ converges
strongly in $ C[0,1] $  to certain $\wt \vphi \in C[0,1]$. This
allows one to pass to limits in (\ref{e int ep})-(\ref{e int bc1})
and to complete the proof using Lemma \ref{l int} .
\end{proof}

\section{Perturbations and derivatives of quasi-eigenvalues.}
\label{s MultEig}

For the sake of convenience, some formal changes in the settings should be done.
We extend the introduced notation to structures $B \in L_\C^\infty (0,1)$.
For complex-valued $B$, the statement of the quasi-eigenvalue problem and
the definition of multiplicities of quasi-eigenvalues remain without changes.

By $\vphi (x, z; B )$ and $\psi (x, z; B)$ we denote the solutions of $y''(x)= - z^2 B(x) y(x) $ satisfying
\begin{equation} \label{e phi psi}
\vphi (0, z; B ) = \pa_x \psi (0, z; B) = 1 , \ \
\pa_x \vphi (0, z; B ) = \psi (0, z; B) = 0 .
\end{equation}

Recall that $F (z;B) := \vphi (1,z;B) - \i \pa_x \vphi (1,z; B)/z$,
and that the set $K_r (B)$ of quasi-eigenvalues of multiplicity $r$
is defined as the set of $r$-fold zeroes of $F(\cdot,B)$.

Basic definitions concerning analytic maps in Banach spaces may be found, e.g., in \cite{PT87}.

\begin{lem} \label{l an}
The map $(z, B ) \mapsto \vphi (\cdot ,z; B)$ is
analytic from $\C \times L_\C^\infty (0,1)$ to
$W_\C^{2,\infty} [0,1]$. Its Maclaurin series is
\begin{eqnarray}
\vphi (x , z ; B ) = 1 - \vphi_1 (x; B ) z^2 + \vphi_2 (x; B ) z^4 - \vphi_3 (x; B )  z^6 + \dots , \label{e pow ser n}\\
\vphi_0 (x; B) \equiv 1 , \ \  \vphi_{j} (x; B ) = \int_{0}^x
(x - s) \, \vphi_{j-1} (s; B ) B (s) \d s , \ \ j \in \N .
\notag 
\end{eqnarray}
\end{lem}

\begin{proof}
The $W_\C^{2,\infty} [0,1]$-valued series (\ref{e pow ser n}) for the solution $\vphi $
is well-known (see e.g. \cite{KK68_II} or \cite{KN89}).
It follows from the estimates in \cite[Sec.2]{KK68_II} (see also \cite[Exercises
5.4.2-3]{DM76}) that the series converge uniformly on every bounded
set of $\C \times L^\infty_\C$. So $(z, B ) \mapsto \vphi (\cdot ,z; B)$ is an analytic map on
$\C \times L^\infty_\C$ (see e.g. \cite[Theorem A.2]{PT87}).
\end{proof}

\begin{lem} \label{l der F}
(i) $F (z;B)$ is  analytic on $\C \times L_\C^\infty
(0,1)$.
\item[(ii)] At quasi-eigenvalues $\k \in K (B)$, the derivative of $ F $ w.r.t. $z$ is given by
\begin{equation} \label{e paz F}
\pa_z F (\k;B) = 2 \left[- \k \psi (1,\k;B) + \i \pa_x \psi (1,\k;B)  \right] \int_0^1 \vphi^2 (s,\k;B) B (s) \d s
+ \frac{\vphi (1,\k;B)}{\k}
;
\end{equation}
and the directional derivatives $[\pa_{B} F (\k,B)] \ (B_\Delta)$ w.r.t. $B$ in the direction $B_\Delta \in L_\C^\infty (0,1)$ by
\begin{equation} \label{e paBD F}
[\pa_{B} F (\k,B)] (B_\Delta) = \k \left[- \k \psi (1,\k;B) + \i \pa_x \psi (1,\k;B) \right] \int_0^1
\vphi^2 (s,\k;B) \ B_\Delta (s) \ \d s .
\end{equation}
\end{lem}

\begin{proof}

\textbf{(i)} follows from Lemma \ref{l an}.

\textbf{(ii)} \emph{Differentiation w.r.t. $B$.}
For any $z \in \C$, the functions
$\vphi (x,z;B)$ and $\pa_x \vphi (x,z;B)$ satisfy
\begin{eqnarray}
\vphi (x,z;B) & = & 1 - z^2 \int_0^x (x-s) \ B (s) \ \vphi (s,z;B) \ \d s , \label{e phi z int}\\
\pa_x \vphi (x,z;B) & = & -z^2 \int_0^x \ B(s) \ \vphi (s,z;B) \ \d s . \label{e pax phi z int}
\end{eqnarray}
To find directional derivatives $[\pa_{B} \vphi] (B_\Delta)$ and $[\pa_{B} \pa_x \vphi] (B_\Delta)$,
we differentiate  these equalities by definition using Lemma \ref{l an}.
We get
\begin{gather}
[\pa_{B} \vphi (x,z;B)] (B_\Delta) = -
z^2 \int_0^x (x-s) B (s) [\pa_{B} \vphi (x,z;B)]
( B_\Delta ) \d s 
- z^2 \int_0^x (x-s) B_\Delta (s)  \vphi (x,z;B) \d s,
  \label{e int paBD ph} \\
 \left[ \pa_{B} \pa_x \vphi (x,z;B) \right] (B_\Delta)  =  -z^2 \int_0^x B(s) [\pa_{B} \vphi (s,z;B)] (B_\Delta) \d s -z^2 \int_0^x B_\Delta (s) \vphi (s,z;B) \d s.
 \label{e int paBD pax ph}
\end{gather}
It follows from (\ref{e int paBD ph}) that $[\pa_{B} \vphi (x,z;B)] (B_\Delta)$
is the solution $y(x)$ of the initial value problem
\begin{equation} \label{e nh bvp}
 y''(x) + z^2 B(x) y(x) = f(x) ,  \ \ y(0) = 0, \ y'(0) = 0
\end{equation}
with $ f(x) = - z^2 B_\Delta (x) \vphi (x,z;B)$.
Further, (\ref{e int paBD pax ph}) can be rewritten as
$[\pa_{B} \pa_x \vphi (x,z;B)] (B_\Delta) = y'(x)$
(as a by-product, we get
$ [\pa_{B} \pa_x \vphi] (B_\Delta) = \pa_x [\pa_{B} \vphi ] (B_\Delta)$).
Solving (\ref{e nh bvp}) by variation of parameters, one can find $y$,$y'$, and, in turn, $y(x) - \i y'(x)/z$.
For $z \neq 0$, $y(x) - \i y'(x)/z$ equals
\begin{multline*} 
\int_0^x f(s) \Bigl( \bigl[ \psi (x,z;B) - \i \pa_x \psi (x,z;B) /z \bigr] \vphi (s,z;B) - \bigl[ \vphi (x,z;B) - \i \pa_x \vphi (x,z;B) /z \bigr] \psi (s,z;B) \Bigr) \d s .
\end{multline*}
Substituting $x=1$, $f$, and $z=\k \in K (B)$ (so that $\vphi (1,\k;B) - \i \pa_x \vphi (1,\k;B) /\k = 0$), we get (\ref{e paBD F}).

\emph{Differentiating (\ref{e phi z int}), (\ref{e pax phi z int}) w.r.t. }$z$,
we see that $\pa_z \vphi $
is given by the solution $y$ of (\ref{e nh bvp}) with $ f(x) = - 2 z B (x) \vphi (x,z;B)$ and that $\pa_z \pa_x \vphi = \pa_x \pa_z \vphi $. Hence, for $z=\k \in K (B)$,
\begin{eqnarray*}
\pa_z F (\k;B) & = & \pa_z \vphi (1,\k;B) - \frac{\i \pa_z \pa_x \vphi (1,\k;B)}{\k} + \frac{\i \pa_x \vphi (1,\k;B)}{\k^2} = y(1) - \frac{\i y' (1)}{\k} + \frac{\i \pa_x \vphi (1,\k;B)}{\k^2}
\\ & = &
 2 \left[- \k \psi (1,\k;B) + \i \pa_x \psi (1,\k;B)  \right] \int_0^1 \vphi^2 (s,\k;B) B (s) \d s + \frac{\i \pa_x \vphi (1,\k;B)}{\k^2}.
\end{eqnarray*}
Using
$ \vphi (1,\k;B) - \i \pa_x \vphi (1,\k;B) /\k =0 $ to modify the last term, we get (\ref{e paz F}).
\end{proof}

Since the solutions $\vphi (\cdot,z;B)$ and $\psi (\cdot,z;B)$ are linearly independent,
\begin{equation} \label{e phi psi alt}
\text{at most one of the numbers } F(z;B) \text{ and } \left[- z \psi (1,z;B) + \i \pa_x \psi (1,z;B)  \right] \ \text{can be } 0 .
\end{equation}
In particular, $\pa_{B} F (\k,B) \neq 0$ if $\k \in K(B)$.

\begin{prop}[cf. \cite{HBKW08} and the discussion in Sec.5 of \cite{AAD01}] \label{p paB k}
Let $\k_0 \in K_1 (B_0)$ (i.e., $\k_0$ is a simple quasi-eigenvalue). Then there exist an open ball $\B_{ \rho} (B_0) \subset L^\infty_\C (0,1)$, $\rho>0$, and a unique continuous functional $k : \B_{\rho} (B_0) \to \C$, such that $k (B) \in K (B)$ and $k (B_0) = \k_0$. Moreover, $k$ is analytic in $\B_{\rho} (B_0)$ and
\begin{equation} \label{e der k}
[\pa_{B} k (B)] (B_\Delta)
= - \frac{\k_0^2 \int_0^1
\vphi^2 (s,\k_0;B) \ B_\Delta (s) \ \d s}{2 \k_0 \int_0^1 \vphi^2 (s,\k_0;B) B (s) \d s -  \i \vphi^2 (1,\k_0;B)} .
\end{equation}
\end{prop}

\begin{proof}
$\k_0 \in K_1 (B_0)$ if and only if $F(\k_0;B_0) = 0$ and
$\pa_z F (\k_0,B_0) \neq 0$. The existence of the functional $\k_0 (B)$ with the desired properties follows from Lemma \ref{l der F} (i) and the implicit function theorem for analytic maps (see e.g. \cite[Appendix B]{PT87}). To get
(\ref{e der k}), we differentiate $F(k(B);B) = 0$ w.r.t. $B$ in the direction $B_\Delta$, and then, use Lemma \ref{l der F} (ii) and the equalities
$\left| \begin{array}{cc}
\vphi (x) & \psi (x) \\
\vphi' (x) & \psi' (x)
\end{array} \right| = 1$, \  $\vphi (1,\k_0;B) = \i \pa_x \vphi (1,\k_0;B) /\k_0 $.
\end{proof}

\begin{lem} \label{l P ser}
Let $P(z,\zeta) = z^r + h_1 (\zeta) z^{r-1} + \dots + h_r (\zeta)$ be a monic polynomial in $z$ with coefficients
$h_j$ analytic in $\zeta$ in a neighborhood of $\zeta_0 = 0$. Suppose $h_j (0) = 0$, $j = 1, \dots, r$, and $h_r '(0) \neq 0$. Then for $\zeta$ close to $0$ there exist
exactly $r$ distinct roots of $P (\cdot,\zeta)= 0$ and these roots are given by an
 $r$-valued analytic  function $Z (\zeta)$ that admits a Puiseux series representation
\begin{equation} \label{e Z Pser}
Z (\zeta) = \sum_{j=1}^{\infty} c_j \zeta^{j/r} \  \
\text{ with } \ c_1 = \sqrt[\r]{-h'_r (0)} \neq 0
\end{equation}
($\sqrt[\r]{\zeta}$ is an arbitrary fixed branch of the multi-function $\zeta^{1/r}$, $c_j \in \C$ are constants).
\end{lem}

\begin{proof}
By the implicit function theorem for multiple zeroes (see e.g. \cite[Theorem XII.2]{RS78IV}), there exist a natural number $p\leq r$ and a (possibly multi-valued) convergent Puiseux series
$ Z (\zeta) = \sum_{j=1}^{\infty} c_j \zeta^{j/p} $ such that all its values are roots of
$P (\cdot,\zeta)= 0$ for $\zeta$ small enough.

Since $Z (\zeta) = o (1)$ as $\zeta \to 0$ and $h_j (0) = 0$, we see that
\[
Z^r (\zeta) = - h_r (\zeta) - h_{r-1} (\zeta) Z (\zeta) - \dots - h_{1} (\zeta) Z^{n-1} (\zeta)
= - \zeta [h'_r(0) + o(1)]  .
\]
Therefore $Z (\zeta) = [-\zeta h_r '(0)]^{1/r} + o (|\zeta|^{1/r})$ with $h'_r(0) \neq 0$ for at least one branch of
$\zeta^{1/r} $ (the equality is valid in domains with a cut, e.g., $\D_\delta (0) \setminus e^{\i \xi }  \R_+ $ ).
This implies $p \geq r$, and so, $p=r$.
Further, the first term $ c_1 \zeta^{1/r}$ is $[-\zeta h_r '(0)]^{1/r} \not \equiv 0$.
So $Z$ gives exactly $r$ distinct roots of $P (\cdot,\zeta)= 0$.
Since $P (\cdot,\zeta)= 0$ has at most $r$ roots, we obtain the statement of the lemma.
\end{proof}

The following proposition describes splitting of an $r$-fold quasi-eigenvalue $\k_0$ under perturbations $\zeta B_\Delta$ satisfying $[\pa_B F (\k_0; B_0)](B_\Delta) \neq 0$. It is essential that in this case the multiple quasi-eigenvalue splits into simple quasi-eigenvalues
'uniformly' in $r$ directions (like roots of $z^r + C\zeta=0$, $C \neq 0$).

\begin{prop} \label{p per k mult}
Let $\k_0 \in K_r (B_0)$ with $r \geq 2$. Assume that
$B_\Delta \in L_\C^\infty$ satisfies
\begin{equation} \label{e paB neq 0}
\int_0^1 \vphi^2 (s,\k_0;B_0) \ B_\Delta (s) \ \d s \neq 0.
\end{equation}
 Then there exist open discs $\D_{\delta} (0) , \D_{\ep} (\k_0) \subset \C $, $\delta,\ep>0$,
 and a convergent in $\D_{\delta} (0)$ $r$-valued Puiseux series
\begin{equation} \label{e k Pser}
k (\zeta) = \k_0 + \sum_{n=1}^{\infty} c_j \zeta^{j/r} \
\text{ with } c_1 = \sqrt[\r]{-\frac{r! \, [\pa_B F ( \k_0 ; B_0 )] (B_\Delta)}{\pa_z^r F (\k_0; B_0)}} \neq 0
\end{equation}
such that for any $\zeta \in \D_{\delta} (0)$, the $r$ values of $k
(\zeta)$ give all the quasi-eigenvalues of $B_0 + \zeta B_\Delta$ in
$\D_{\ep} (\k_0)$ and all these $r$ quasi-eigenvalues are distinct
and simple.
\end{prop}

\begin{proof}
Consider the entire function $\wt F (z,\zeta) := F ( z ; B_0 + \zeta
B_\Delta )$ of two complex variables $z$ and $\zeta$. Then $\k_0$ is
an $r$-fold zero of the function $\wt F (\cdot,0)$. By the
Weierstrass preparation theorem, in a certain polydisc $\D_{\ep_1}
(\k_0) \times \D_{\delta_1} (0)$,
\[
\wt F (z,\zeta) = 
[(z-\k_0)^r + h_1 (\zeta) (z-\k_0)^{r-1} + \dots + h_r (\zeta) ] G (z,\zeta) ,
\]
where the coefficients $h_j$ (the function $G (z,\zeta)$) are
analytic in $\D_{\delta_1} (\k_0)$ (resp., in $\D_{\ep_1} (\k_0)
\times \D_{\delta_1} (0)$), $h_j (0) = 0$,  and $G (z,\zeta) \neq 0
$ in $\D_{\ep_1} (\k_0) \times \D_{\delta_1} (0)$.

It follows from (\ref{e phi psi alt}), (\ref{e paB neq 0}), and Lemma \ref{l der F} (ii) that $[\pa_B F ( \k_0 ; B_0 )] (B_\Delta) \neq 0$.
Differentiating $\wt F$ by $\zeta$, one gets
\[
[\pa_B F ( \k_0 ; B_0 )] (B_\Delta) = \pa_\zeta \wt F (\k_0,0)  = h'_r (0) G (\k_0,0)  .
\]
On the other side, $\pa_z^r \wt F (\k_0,0) = r! G (\k_0,0) $. Hence,
\[
h'_r (0) = \frac{[\pa_B F ( \k_0 ; B_0 )] (B_\Delta)}{ G (\k_0,0)} = \frac{r! \, [\pa_B F ( \k_0 ; B_0 )] (B_\Delta)}{\pa_z^r F (\k_0,B_0)} \neq 0.
\]
So Lemma \ref{l P ser} may be applied to the Weierstrass polynomial
$ (z-\k_0)^{r} + h_1 (\zeta) (z-\k_0)^{r-1} + \dots + h_r (\zeta)$
to get the Puiseux series (\ref{e k Pser}) for its
zeroes, which are also zeroes of
$\wt F$, and so, are quasi-eigenvalues of $B_0 + \zeta B_\Delta $. In this way, we get $r$ distinct
quasi-eigenvalues of $B_0 + \zeta B_\Delta $ that approach $\k_0$ as $\zeta \to 0$.

Since the obtained $r$ zeroes of $\wt F (\cdot,\zeta)$ are distinct for small $\zeta$,
the standard Rouche's theorem argument implies that each of them is of multiplicity 1.
\end{proof}

The study of two-parameter perturbations of quasi-eigenvalues requires the following lemma.

Let us denote
\begin{eqnarray*}
T_0 & := & \{\zeta = (\zeta_1,\zeta_2) \in \R^2\, : \, \zeta_1, \zeta_2 \geq 0
\ \text{ and } \ \zeta_1 + \zeta_2 \leq 1 \} . \label{e T0}
\end{eqnarray*}

\begin{lem} \label{l 2par per}
Let $Q (z,\zeta_1,\zeta_2) $ be a function of three complex variables analytic in a neighborhood of the origin
$\mathbf{0}=(0,0,0)$. Assume that $0$ is an $r$-fold zero ($1\leq r < \infty$) of the function
$Q (\cdot,0,0)$, that $\pa_{\zeta_{\scriptstyle j}} Q (\mathbf{0}) \neq 0$, $j=1,2$, and
$ \arg \pa_{\zeta_{\scriptstyle 2}} Q (\mathbf{0}) =  \arg \pa_{\zeta_{\scriptstyle 1}} Q (\mathbf{0}) + \xi_0 \ (\modn 2\pi)$ with
$\xi_0 \in (0,\pi)$. Denote $\eta_j := -\frac{r! \pa_{\zeta_{\scriptstyle j}} Q (\mathbf{0})}{\pa_z^r Q (\mathbf{0})}$, $j=1,2$.
Then for any $\delta>0$ and $\xi_1$ in the interval
$\left( \arg \sqrt[r]{\eta_1} , \frac{\xi_0}{r} + \arg \sqrt[r]{\eta_1} \right)$, there exist a pair
$(\zeta_1,\zeta_2) \in T_0 \cap (\D_\delta (0) \times \D_\delta (0))$ and
$z \in \C \setminus \{0\}$ such that $Q (z,\zeta_1,\zeta_2) = 0$ and\\
 $\arg z = \xi_1 \ (\modn 2\pi)$.
\end{lem}

\begin{proof}
We give the detailed proof for the more difficult case
$r \ge 2$. The case when $0$ is a simple zero of $Q (\cdot,0,0)$ is similar, but much simpler in notation and details.

It is not an essential restriction to assume that $\sqrt[\r]{\cdot}$ is the analytic in $\C \setminus \overline{\R_-} $ branch of $(\cdot)^{1/r}$ fixed by  $\sqrt[\r]{1}=1$, and that
\begin{eqnarray}
& \arg \pa_{\zeta_j} Q (\mathbf{0}) = (-1)^j \frac{\xi_0}{2} , \ \ \ \arg\sqrt[r]{\eta_j} = (-1)^j \frac{\xi_0}{2r}, \ \ j=1,2, \label{e arg pa Q}
\end{eqnarray}
($\pa_z^r Q (\mathbf{0}) $ can be placed on $\R_-$ by a change of variable in $z$).

In this settings, we have to prove that for any
$\xi_1 \in \left( - \frac{\xi_0}{2r} , \frac{\xi_0}{2r} \right)$
there exists a sequence $\{ (z_n , \zeta_1^{[n]} , \zeta_2^{[n]}) \}$ going to $\mathbf{0}$ and such that $\arg z_n = \xi_1$, $\{ \zeta^{[n]} \} \subset T_0$, and $ Q(z_n, \zeta^{[n]}) = 0 $.

\emph{Step 1. An auxiliary triangle and the Weierstrass
decomposition.} Let us introduce the triangle
\begin{eqnarray*}
 T_1 & = & \left\{ \left([1-\theta] c , \theta c\right) \in \R^2 \, : \, c \in [0,1], \ 0< \theta_1 \leq \theta \leq \theta_2 < 1 \right\} ,
\end{eqnarray*}
with $\theta_1, \theta_2$ such that
\begin{eqnarray} \label{e <xi1<}
\arg \sqrt[\r]{\wt \eta_1} & < & \xi_1 < \arg \sqrt[\r]{\wt \eta_2},
\ \text{ where } \wt \eta_j := (1-\theta_j) \eta_1 + \theta_j \eta_2 , \ \
j=1,2.
\end{eqnarray}
 Clearly, $T_1 \subsetneqq T_0$, and
if a sequence $\{ \zeta^{[n]} \} \subset T_1 $ tends to $ (0,0)$, we have
\begin{equation} \label{e ga asymp}
\zeta_1^{[n]} \asymp \zeta_2^{[n]} \asymp |\zeta^{[n]}| \ \text{ as } n \to \infty .
\end{equation}

By the Weierstrass preparation theorem, in a certain polydisc
$\D_{\ep_1} (0) \times \D_{\delta} (0) \times \D_{\delta} (0)$,
\begin{equation} \label{e Qzga}
Q (z,\zeta) = P (z,\zeta) R (z,\zeta) \
 \ \text{ with } P (z,\zeta) = z^r + q_1 (\zeta) z^{r-1} + \dots + q_r (\zeta) ,
\end{equation}
where the coefficients $q_j$ of the Weierstrass polynomial $P$ (the
function $R$) are analytic in $\D_{\delta} (0) \times \D_{\delta}
(0)$ (resp., in $\D_{\ep_1} (0) \times \D_{\delta} (0) \times
\D_{\delta} (0)$), $q_j (0,0) = 0$, and $R (z,\zeta_1,\zeta_2) \neq
0 $ in $\D_{\ep_1} (0) \times \D_{\delta} (0)  \times \D_{\delta}
(0)$. We can suppose that $\delta >1$. Indeed, we can always ensure
this scaling the variables $\zeta_1$, $\zeta_2$. Now the Weierstrass
decomposition (\ref{e Qzga}) holds for $(\zeta_1,\zeta_2) \in T_0$.

\emph{Step 2. Asymptotics of zeroes of $Q(z,\zeta)$ for $\zeta \in
T_1$.}

Using the arguments of Proposition \ref{p per k mult},
one can show that
\begin{equation} \label{e pa qr}
\pa_{\zeta_j} q_r (0,0) = \frac{\pa_{\zeta_j} Q (\textbf{0}) }{R (\textbf{0})} = \frac{r! \, \pa_{\zeta_j} Q (\textbf{0})}{\pa_z^r Q (\textbf{0})} \neq 0 , \ \ j=1,2.
\end{equation}

Assume that there exist $z_n \to 0$ and $(\zeta_1^{[n]},\zeta_2^{[n]})  \to (0,0)$ as $n \to \infty$
such that $Q (z_n,\zeta^{[n]}) = 0$ and $\{ \zeta^{[n]} \} \subset T_1$. Then it follows from (\ref{e Qzga}) that
\[
z_n^r = - q_r (\zeta^{[n]}) + q_{r-1} (\zeta^{[n]}) o (1) + \dots + q_{1} (\zeta^{[n]}) o (1) , \ \ n \to \infty .
\]
Using (\ref{e pa qr}), (\ref{e ga asymp}), and $q_j (0,0) = 0$, one can show that
\begin{equation} \label{e as zn}
z_n = \left(- \sum_{j=1}^2 \zeta_j^{[n]} \pa_{\zeta_j} q_r (0,0) \right)^{1/r}  [1 + o (1)]
\asymp |\zeta^{[n]}|^{1/r} , \ \ n \to \infty .
\end{equation}

\emph{Step 3. Multiplicities of zeroes of $Q(z,\zeta)$ for $\zeta
\in T_1$.} Let us show that for $\zeta \in T_1$ small enough and $z$ small enough, roots
of $Q(z,\zeta) = 0$ are simple.

Assume the contrary. Then there exist sequences $z_n $ and
$(\zeta_1^{[n]};\zeta_2^{[n]}) $ as above with the additional
property that $\pa_z Q(z_n,\zeta^{[n]}) = 0$. This implies $\pa_z P
(z_n,\zeta^{[n]}) = 0$. Using (\ref{e ga asymp}) and $q_j (0,0) = 0$
again, we see that
\[
r z_n^{r-1} = q_{r-1} (\zeta^{[n]}) + \dots + (r-1) q_{1} (\zeta^{[n]}) z_n^{r-2} = |\zeta^{[n]}| \, O (1), \ \ n \to \infty .
\]
So $z_n = |\zeta^{[n]}|^{1/(r-1)} O (1) $. This contradicts
(\ref{e as zn}).

Rescaling $\zeta_1,\zeta_2$ if necessary, we can ensure that all the roots of $Q(z,\zeta) = 0$
are simple for $\zeta \in T_1$ .

\emph{Step 4.} Applying arguments of Proposition \ref{p per k mult}
 to the zeroes of the function $Q(\cdot , (1-\theta) \tau , \theta \tau)$ with  a complex variable $\tau$
 and a fixed parameter $\theta \in [\theta_1,\theta_2]$ (for example, $\theta= \theta_1$), one can
produce the $r$-valued Puiseux series
\begin{equation} \label{e Ztau}
Z ([1-\theta] \tau , \theta \tau) = \sum_{j=1}^{\infty} c_j \tau^{j/r} \
\text{ with } c_1 = \sqrt[\r]{ (1-\theta) \eta_1 + \theta \eta_2 } \neq 0,
\end{equation}
 ($\theta_j$
and $\eta_j$ are from Step 1 and the statement of the lemma). One of
the values $Z_1 ([1-\theta] \tau , \theta \tau) $ of this function
can be chosen if we place the branch $ \sqrt[r]{\tau} $ instead of
the multi-function $\tau^{1/r}$ in the Puiseux series (\ref{e
Ztau}). Step 3 and the implicit function theorem for simple zeros
imply that $Z_1$ can be extended to an analytic on $T_1\setminus \{ 0 \}$ function
$Z_1 (\zeta_1,\zeta_2)$. For each $\theta \in [\theta_1,\theta_2]$
and sufficiently small $\tau \in (0, \epsilon_1 (\theta)]$,
$\epsilon_1 (\theta)>0$, the function $Z_1 ([1-\theta] \tau , \theta
\tau) $ is given by the Puiseux series (\ref{e Ztau}) with
$\tau^{1/r}$ replaced by its branch $\sqrt[\r]{\tau}$ and
coefficients $c_j$ depending on $\theta$. (However, we do not know
if $\epsilon (\theta)$ is uniformly separated from $0$, and so, we
use other arguments to study the asymptotics of $Z_1$).

Comparing (\ref{e Ztau}) with (\ref{e as zn}) and using (\ref{e ga asymp})
, one can see that
\begin{equation} \label{e z1}
Z_1 (\zeta_1 , \zeta_2) =
\sqrt[\r]{ \zeta_1 \eta_1 + \zeta_2 \eta_2 } [1 + o(1)] , \ \ \zeta \to 0, \ \ \zeta \in T_1 .
\end{equation}
This implies that for $\zeta \in T_1$ small enough,
$Z_1 (\zeta) \neq 0$ and, due to (\ref{e arg pa Q}),
$\arg Z_1 (\zeta)$ is a continuous function with values in $(-\pi/2 , \pi/2)$. On the other side, (\ref{e as zn}) and (\ref{e <xi1<}) imply that for $\tau$ small enough,
$\arg Z_1 ([1-\theta_1] \tau , \theta_1 \tau) < \xi_1 $ and $\arg Z_1 ([1-\theta_2] \tau , \theta_2 \tau) > \xi_1 $.
So for each $\tau \in (0,\epsilon_2)$, $\epsilon_2>0$, there exists
$\theta \in (\theta_1,\theta_2) $ such that $\arg Z_1 ([1-\theta] \tau , \theta \tau) = \xi_1 $. This completes the proof.
\end{proof}

\section{Proofs of Theorems \ref{t B in Exts} and \ref{t J alpha0}.}
\label{s pr t extr p}

\subsection{Optimal structures are extreme points of $\Ad$; the case $\alpha \neq 0$, $b_1>0$. \label{ss case aneq0}}

We start from the case when $b_1 > 0$. Note that this assumption is
satisfied in the optimization problem for photonic crystals
\cite{KS08,HBKW08}. For the degenerate case $b_1 =0$, some details
of the proof require modifications (see Subsection \ref{ss case
b1=0} below).

\begin{lem} \label{l phi values}
Let $B \in L_\R^1 (0,1)$, $B(x) > 0 $ a.e., and $z^2 \not \in \R$. Then:
\item[(i)] $\vphi (x,z;B) \neq 0$ for all $x \in [0,1]$.
\item[(ii)] $\frac{\pa_x \vphi (x,z;B)}{\vphi (x,z;B)} \not \in \R $ for all $x \in (0,1]$.
\item[(iii)] For any $\xi \in (-\pi,\pi]$
the set $\{ x \in [0,1] \, : \, \vphi (x,z;B) \in e^{\i \xi} \R_+ \}$ is finite.
\end{lem}

\begin{proof}
\textbf{(i)} Let $\vphi (x_1,z;B) = 0$. Then $x_1 > 0$ (since $\vphi (0,z;B) =1 $), and so, $\vphi (x,z;B)$ is an eigenfunction of the self-adjoint boundary value problem $ -y'' = \la B y$, $y'(0)=y(x_1)= 0$.
Hence the corresponding eigenvalue $\la = z^2 $ is real, a contradiction.

\textbf{(ii)} Let $h = \frac{\pa_x \vphi (x_1,z;B)}{\vphi (x_1,z;B)} \in \R $ for $x_1>0$. Then $\vphi (x,z;B)$ is an eigenfunction of the self-adjoint problem
$ -y'' = \la B y$, $y'(0)=0$, $y'(x_1) - h y(x_1) =0$.
So $ z^2 \in \R$, a contradiction.

\textbf{(iii)} Let $\Omega := \{ x \in [0,1] \, : \, \vphi (x,z;B)
\in e^{\i \xi} \R_+ \}$ be infinite. Then it includes a convergent
sequence of distinct points $\{ x_n \}_{n \in \N} \subset \Omega$.
It follows from statement (i) and $\vphi \in C [0,1]$ that $x_0 :=
\lim x_n$ also belongs to $\Omega$. Then $\frac{\vphi (x_n) - \vphi
(x_0)}{x_n - x_0} \in e^{\i \xi} \R$, and so, $\pa_x \vphi (x_0) \in
e^{\i \xi} \R$. If $x_0 > 0$, this \emph{contradicts statement
(ii)}.

\emph{So $x_0 = 0$.} Assume $z^2 \in \C_-$ (the case $z^2 \in \C_+$
is similar). Then there exists a neighborhood of $(-z^2)$ lying in a
sector separated from $\R$, i.e., there exist an interval $[\xi_1,
\xi_2] \subset (0,\pi)$ and $\ep >0$ such that
\[
\overline{\D_\ep (-z^2)} \subset
\{ \zeta \in \C\setminus{0} \, : \, \arg \zeta \in [\xi_1, \xi_2] \} =: \Sec (\xi_1,\xi_2)  .
\]
By continuity of $\vphi (x) = \vphi (x,z;B)$, there exists $\delta >0$ such that $\{ -z^2 \vphi (x) \, : x \in [0,\delta] \} \subset \overline{\D_\ep (-z^2)}$. Since $B>0$ a.e., we have
$-z^2 B (x) \vphi (x) \in \Sec (\xi_1,\xi_2)$ for a.a.
$x \in [0,\delta]$.
Using (\ref{e phi z int}), one can show that $\vphi (x) \in 1+  \Sec (\xi_1,\xi_2)$ for all $x \in (0,\delta]$.
The intersection of the ray $e^{\i \xi} \R_+$ with $1+ \Sec (\xi_1,\xi_2)$ is either empty or
separated from the point $1$ (note that $0 \not \in \Sec (\xi_1,\xi_2) $).
Using $\vphi(0)=1$ and $\vphi \in C [0,1]$, we can choose $\delta$ so small that
$ \vphi (x) \not \in e^{\i \xi} \R_+$ for all $x \in (0,\delta]$. So $x_0=0$ is not a limit point of $\Omega$, a contradiction.

\end{proof}

\begin{lem} \label{l paB Ad conv}
Let $B_0 \in \Ad $, $z \in \C$. Then the set of directional derivatives
$
 [\pa_B F (z;B_0)] (\Ad - B_0) := \{ [\pa_B F (z;B_0)] (B_\Delta) \, : \, B_\Delta + B_0 \in \Ad \}
$
is convex.
\end{lem}

The lemma follows from the fact that $\Ad - B_0 := \{ B_\Delta : B_\Delta + B_0 \in \Ad \}$ is convex.

\begin{lem} \label{l 0 int}
Let $B_0 \in \Ad \setminus \Ext $, $\k_0 \in K (B_0)$, and $\re k_0 \neq 0$. Then $0$ is an interior point of the set $[\pa_B F (\k_0;B_0)] (\Ad - B_0) $.
\end{lem}

\begin{proof}[Proof for the case $b_1 >0$.]
Since $B_0 \in \Ad \setminus \Ext $, there exist $\ep_1 >0$ such that the set
 $\Omega := \{ x \in (0,1) \, : \, b_1+\ep_1 < B_0 (x) < b_2 -\ep_1\}$ is of positive Lebesgue measure.
 Using Lemma \ref{l phi values}, one can show that there exist two sectors
$ \Sec (\xi_1,\xi_2) $ and $\Sec (\xi_3,\xi_4) $ such that
the sets
\[
\Omega_1 := \{ x \in \Omega \, : \, \vphi^2 (x,\k_0;B_0) \in \Sec (\xi_1,\xi_2) \}  \text{ and }
 \ \Omega_2 := \{ x \in \Omega \, : \, \vphi^2 (x,\k_0;B_0) \in \Sec (\xi_3,\xi_4) \}
\]
are of positive  measure, and
\[
\Sec (\xi_1,\xi_2) \cap \Sec (\xi_3,\xi_4) = \emptyset , \ \
\Sec (\xi_1,\xi_2) \cap (-1) \Sec (\xi_3,\xi_4) = \emptyset .
\]
This, $\k_0 \in K (B_0)$, and (\ref{e phi psi alt}) imply that the complex points $\zeta_1$ and $\zeta_2$ defined by
\[
\zeta_j  := \k_0 \left[- \k_0 \psi (1,\k_0;B_0) + \i \pa_x \psi (1,\k_0;B_0) \right] \int_{\Omega_j} \vphi^2 (s,\k_0;B_0) \d s , \ \ j=1,2,
\]
are linearly independent as vectors in $\R^2$. Due to Lemma \ref{l der F} (ii), $\zeta_j = [\pa_{B} F (\k_0,B_0)] (\chi_{\Omega_j})$ .

Taking $B_{\Delta,j} (x) = \ep\chi_{\Omega_j} (x) $, $\ep \in (0,\ep_1)$, $j=1,2$, we see that $B_{\Delta,j} \in \Ad - B_0 $ (by the definition of the set $\Omega$). So the complex intervals $\zeta_j (-\ep_1,\ep_1)$, $j=1,2$, are subsets of $[\pa_B F (\k_0;B_0)] (\Ad - B_0) $.

Since $\zeta_1$ and $\zeta_2$ are linearly independent, $0$ is an interior point of the convex hull of the
intervals $\zeta_1 (-\ep_1,\ep_1) $ and $\zeta_2 (-\ep_1,\ep_1) $. Lemma \ref{l paB Ad conv} concludes the proof.
\end{proof}

\begin{prop} \label{p k nonmin}
Let $B_0 \in \Ad $ and $\k_0 \in K (B_0)$.
If $0$ is an interior point of the set\\ $[\pa_B F (\k_0;B_0)] (\Ad - B_0) $, then there
exists $\beta_1>0$ such that $\k_0 - \i \beta_1  \in K (\Ad)$.
\end{prop}

\begin{proof}
The quasi-eigenvalue $\k_0$ of $B_0$ has a finite multiplicity $r$ (see e.g. \cite{KN89,CZ95}).
Recall that $\k_0 \in K_r (B_0)$ implies $\pa_z^r F (\k_0; B_0) \neq 0$.

Since $0$ is an interior point of $[\pa_B F (\k_0;B_0)] (\Ad - B_0) $, there exist
$B_{\Delta,1}, B_{\Delta,2} \in \Ad - B_0$
and a branch $\sqrt[\r]{\zeta}$ of multi-
function $\zeta^{1/r}$
such that
\begin{itemize}
\item[(i)] The complex numbers $\displaystyle \eta_j := -\frac{r! \, [\pa_B F ( \k_0 ; B_0 )] (B_{\Delta,j}) }{\pa_z^r F ( \k_0 ; B_0 )}$, $j=1,2$, are nonzero.
\item[(ii)] $\arg \eta_2 = \arg \eta_1 + \frac{\pi}{2} \ (\modn 2 \pi)$.
\item[(iii)] $\sqrt[\r]{\cdot}$ is holomorphic in the sector $\Sec_1 := \{ c_1 \eta_1 + c_2 \eta_2 \, : c_1, c_2  > 0 \} $.
\item[(iv)] $\arg \sqrt[\r]{\eta_j} = -\frac{\pi}{2} + (-1)^j \frac{\pi}{4r} \ (\modn 2\pi)$, $j=1,2$.
\end{itemize}

By Lemma \ref{l der F} (i),  $Q (\wt z,\zeta_1,\zeta_2) := F ( \k_0
+ \wt z; B_0 + \zeta_1 B_{\Delta,1} + \zeta_2 B_{\Delta,2})$ is an
entire function of three complex variables $\wt z$, $\zeta_1$, and
$\zeta_2$. Obviously, $\k_0$ is an $r$-fold zero of $Q (\cdot,0,0)$,
and due to the properties (i)-(iv), we can apply Lemma \ref{l 2par
per} to the function $Q$. Lemma \ref{l 2par per} implies that there
exist $(\zeta_1,\zeta_2) \in T_0 $ and $\wt z \in \C \setminus
\{0\}$ such that $Q (\wt z,\zeta_1,\zeta_2) = 0$ and $\arg \wt z =
-\pi/2 \ (\modn 2\pi)$. So $\wt z = -i\beta_1$, $\beta_1>0$, and
$\k_0 -i\beta_1 \in K ( B_0 + \zeta_1 B_{\Delta,1} + \zeta_2
B_{\Delta,2})$.

Since $ \Ad - B_0$ is convex and contains $0$, the structure $B_0 + \zeta_1 B_{\Delta,1} + \zeta_2 B_{\Delta,2}$
belongs to $\Ad $.
\end{proof}

We see that Lemma \ref{l 0 int} and Proposition \ref{p
k nonmin} imply the desired statement (for the case $b_1>0$).

\begin{cor} \label{c k nonmin}
Let $\alpha \neq 0$, $B \in \Ad $,  and $\alpha + \i \I (\alpha) \in K (B)$.
Then $B \in \Ext$.
\end{cor}


\subsection{Piecewise constancy of optimal structures for $\alpha \neq 0$, $b_1 > 0$. \label{ss case aneq0 step}}

\begin{defn} Let $b_1 \neq b_2$ and $B \in \Ext$.
\item[(i)] We say that $x_0 \in [0,1]$ is a switch point of $B$ if for any $\epsilon>0$
both the sets $E_j (B) \cap (x_0-\epsilon,x_0+\epsilon) $, $j=1,2$, are of positive measure
(see (\ref{e E12}) for the definition of $E_j (B)$).
\item[(ii)] We say that a switch point $x_0$ is singular from the left (right) if
for any $\epsilon>0$
both the sets $E_j (B) \cap (x_0-\epsilon,x_0) $
(resp., $E_j (B) \cap (x_0,x_0+\epsilon) $),
 $j=1,2$, are of positive measure. Otherwise, $x_0$ is called regular from the left (resp., right).
\item[(iii)] A switch point $x_0$ is said to be regular if it is regular both from the left
and from the right.
\end{defn}

It is easy to notice the following.

\begin{lem} \label{l sw seq}
Let $B \in \Ext$ have a switch point $x_0$ singular from the left (right).
Then there exists a sequence $\{ x_n \}_1^\infty$ of (distinct) switch points converging to $x_0$
from the left (resp., right).
\end{lem}

This implies easily the next statement.

\begin{lem} \label{l sw p}
Let $B \in \Ext$. Then the following statements are equivalent:
\item[(i)] $B \in \Exts$.
\item[(ii)] $B$ has a finite number of switch points.
\item[(iii)] All switch points of $B$ are regular.
\end{lem}

Let us define a quasi-argument function $\qarg : \C\setminus\{0\} \to \T $
 by
$
\qarg z = e^{\i \arg z }
$.

Note that for $z \in \C_+ \setminus \i\R$ and $B \in L^1_\R (0,1)$, $B > 0$ a.e.,
Lemma \ref{l phi values} implies that,
for a suitable branch of the argument function,
\begin{equation} \label{e arg vphi2}
 \arg \vphi^2 (x,z;B) \ \text{ is continuously differentiable on } [0,1] \text{ and
} \ \pa_x \arg \vphi^2 (x,z;B) \neq 0 \ \text{ for } \ x\in
(0,1].
\end{equation}
In the next lemma, we keep such a choice of the $\arg$-function.

\begin{prop} \label{p semi-c}
Let $0<b_1$.
Let $B_0 \in \Ext $ have a switch point $x_0 \in (0,1)$.
Let $\k_0 \in K (B_0) \setminus \i\R $,
\[
\xi (x) = \arg \vphi^2 (x,\k_0;B_0), \ \text{ and } \
\xi_1 := \arg \left( \k_0 \left[ - \k_0 \psi (1,\k_0;B_0) + \i \pa_x \psi (1,\k_0;B_0) \right]\right)
\]
($\xi_1$ is well defined due to (\ref{e phi psi alt})).
Then the set of quasi-arguments of directional derivatives
\begin{equation} \label{e Aset}
A(\k_0; B_0):= \qarg \, [\pa_B F (\k_0;B_0)] \Bigl( [\Ad - B_0]
\setminus \{0\} \Bigr)
\end{equation}
contains at least one of the semi-circles $\{ e^{\i [\xi (x_0) +\xi_1+s]} \, : \, s \in (0,\pi)\}$,
$\{ e^{\i [\xi (x_0) + \xi_1 +s]} \, : \, s \in (-\pi,0)\}$.

More precisely, assume that
the measures of the sets $(x_0-\epsilon,x_0) \cap E_{j_{\scriptstyle 1}} (B)$ and
$(x_0,x_0+\epsilon) \cap E_{j_{\scriptstyle 2}}$,
where $j_1,j_2 \in \{ 1, 2 \}$ and $j_1 \neq j_2$, are nonzero for all $\epsilon>0$. Then
\item[(i)]  in the cases \{$\xi'(x_0) < 0$, $j_1=1$, $j_2=2$\} and
\{ $\xi'(x_0) > 0$, $j_1=2$, $j_2=1$ \}, the set
$A(\k_0; B_0) $ contains $\{ e^{\i [\xi (x_0) +\xi_1+s]} \, : \, s \in (0,\pi)\}$,
\item[(ii)] in the cases \{$\xi'(x_0) < 0$, $j_1=2$, $j_2=1$\} and
\{$\xi'(x_0) > 0$, $j_1=1$, $j_2=2$\}, the set
$A(\k_0; B_0) $ contains $\{ e^{\i [\xi (x_0) + \xi_1 +s]} \, : \, s \in (-\pi,0)\}$.
\end{prop}

\begin{proof}
We consider the case when $\xi'(x_0) < 0$ and $j_1=1$, $j_2=2$.
Arguments for all other possible cases are similar.
Then there exist two sequences $\{ L_n \}_1^\infty$, $\{ R_n \}_1^\infty$  of subsets of
$ (0,1)$ such that
\begin{itemize}
\item[(LR1)] All $L_n$ and $R_n$ are of positive measure.
\item[(LR2)] $L_n \subset (x_0 - \frac{1}{n}, x_0 )$, $ R_n \subset (x_0 , x_0 + \frac{1}{n} )$, $n \in \N$.
\item[(LR3)] $B_0 (x)= b_1$ for $x \in L_n$, $B_0 (x)= b_2$ for $x \in R_n$.
\end{itemize}

The assumptions on $L_n$ and $R_n$ imply
that $ (b_2-b_1)\chi_{L_n} $ and $(b_1-b_2)\chi_{R_n}$ are in
$ \Ad - B_0 $.
So the set $ [\pa_B F (\k_0;B_0)] (\Ad - B_0) $ is convex, contains $ \eta_{L,n} :=[ \pa_B F (\k_0;B_0)] \left( [b_2-b_1]\chi_{L_n} \right)$,
$ \eta_{R,n} :=[\pa_B F (\k_0;B_0)] \left( [b_1-b_2]\chi_{ R_n } \right)$,
and, obviously, contains the point $0$.

Now finding from (\ref{e paBD F}) the arguments of $\eta_{L,n}, \eta_{R,n} $ for large $n$,
one can prove the statement of the lemma.
Indeed, for any $\epsilon>0$, taking $n$ large enough, we can ensure that
$\xi (x) \in (\xi (x_0), \xi (x_0)+\epsilon) $ for $x \in L_n$,
and that $\xi (x) \in (\xi (x_0)-\epsilon, \xi (x_0) )$ for $x \in R_n$
(the assumption $\xi'(x_0) < 0$ is used here).
From this and (\ref{e paBD F}) we see that
$
 \arg \eta_{L,n} = \xi_1 + \xi_{L,n}  \ (\modn 2\pi),
$
where
\[
\xi_{L,n} := \arg \left[ (b_2-b_1) \int_{L_n}
\vphi^2 (s,\k_0;B_0) \ \d s \right] \in (\xi (x_0), \xi (x_0)+\epsilon) .
 \]
Similarly, taking into account $b_1 -b_2 <0$,
\[
\arg \eta_{R,n}
= \xi_1 + \xi_{R_n} \ (\modn 2\pi), \ \ \ \
\xi_{R_n}
\in (  \xi (x_0) + \pi  -\epsilon, \xi (x_0) + \pi ) .
 \]
So the quasi-arguments of the convex hull of $0$, $\eta_{L,n} $, and $\eta_{R,n} $
cover at least the set $\{ e^{\i [\xi (x_0) +\xi_1+s ]} \, : \, s \in [\epsilon,\pi-\epsilon]\}$.
Moving $\epsilon \to 0$, we complete the proof.
\end{proof}

\begin{prop} \label{p singsw int0}
Let $B_0 \in \Ext $ have a singular switch point $x_0 \in [0,1]$. Let $\k_0 \in K(B_0)$ and
$\re \k_0 \neq 0$.
Then $0$ is an interior point of the set $ [\pa_B F (\k_0;B_0)] ( \Ad - B_0 ) $.
\end{prop}

\begin{proof}[Proof for the case $b_1>0$.]
By Lemma \ref{l sw seq}, there exists a sequence of distinct switch points $\{ x_j \}_1^\infty$ converging
to $x_0$. Let $\xi (x) = \arg \vphi^2 (x,\k_0;B_0) $ as in Proposition \ref{p semi-c}.
By Lemma \ref{l phi values} (iii), there exist $j$ and $n$ such that
$e^{\i \xi (x_j)} \neq e^{\i \xi (x_n) } $, $e^{\i \xi (x_j)} \neq -e^{\i \xi (x_n) } $, and $x_j, x_n \in (0,1)$.
 Applying Proposition \ref{p semi-c} to the switch points $x_j$ and $x_n$, we see that
$A(\k_0; B_0)$ contains two distinct semi-circles. Since $ [\pa_B F (\k_0;B_0)] (\Ad - B_0) $
is convex and contains $0$, this implies the proposition.
\end{proof}

\begin{cor} \label{c B in Exts}
Let $\alpha \neq 0$, $B_0 \in \Ad $,  and $\k_0 =\alpha + \i \I (\alpha) \in K (B_0)$.
Then $B_0 \in \Exts$.
\end{cor}

\begin{proof}[Proof for the case $b_1 >0$.]
By Corollary \ref{c k nonmin}, $B_0 \in \Ext$. Assume that $B_0 \not \in \Exts$.
Then, by Lemma \ref{l sw p}, $B_0$ has a singular switch point.
Propositions \ref{p singsw int0} and \ref{p k nonmin} imply that there exists
$\k_1 \in K (\Ad )$ such that $\re \k_1 = \alpha $, but
$\im \k_1 < \I (\alpha) $. This contradicts the definition of $\I (\alpha) $.
\end{proof}

\subsection{The case $b_1 = 0$. \label{ss case b1=0}}

In this subsection we complete the proof of Corollary \ref{c B in Exts}
considering the case $b_1 = 0$.

While the most part of the proof remains the same, some changes have
to be done in Lemma \ref{l phi values} and Proposition \ref{p
semi-c}. The case $b_1 = 0$ is degenerate in the sense that the set
$E_1 (B)$ (where $B$ equals $0$) may have a positive measure and
include an interval $[0,x_0)$ with $x_0>0$.

Define
\begin{equation} \label{e a1}
 a_1:= \sup \{ x \in [0,1] \ : \ B = 0 \text{ a.e. on } [0,a_1]  \} .
\end{equation}
So $B(x) =0$ a.e. on $[0,a_1]$, and (if $a_1 <1$)
\begin{equation} \label{e Bneq0 a1x0}
B \ \text{ is not equivalent to } \ 0 \ \text{ on any interval } \ (a_1, x_0) \ \text{ with } \ x_0> a_1.
\end{equation}

Statements (ii) and
(iii) of Lemma \ref{l phi values} are not valid on $[0,a_1]$.
Lemma \ref{l phi values} can be adjusted in the following way.

\begin{lem} \label{l phi values b1=0}
Let $B \in L_\R^1 (0,1)$, $B (x) \geq 0 $ a.e., and $a_1 < 1$. Let $z^2 \not \in \R$. Then:
\item[(i)] $\vphi (x,z;B) \neq 0$ for all $x \in [0,1]$.
\item[(ii)] $\frac{\pa_x \vphi (x,z;B)}{\vphi (x,z;B)} \not \in \R $ for all $x \in (a_1,1]$.
\item[(iii)] For any $\xi \in (-\pi,\pi]$
the set $\{ x \in [a_1,1] \, : \, \vphi (x,z;B) \in e^{\i \xi} \R_+ \}$ is finite.
\end{lem}

\begin{proof}
\textbf{(i)} Obviously, $\vphi (x,z;B) = 1$ for $x \in [0,a_1]$.
Let $\vphi (x_1,z;B) = 0$ for $x_1 > a_1$. Then $\vphi (x,z;B)$
is an eigenfunction of the nontrivial self-adjoint problem defined on
the interval $[a_1,x_1]$ by (\ref{e ep}) and the boundary
conditions $y'(a_1) = y(x_1)= 0$ (see e.g. \cite{KK68_II,DM76}).
The corresponding eigenvalue $\la = z^2 $ is real, a contradiction. In the same way we get statement
\textbf{(ii)}.

\textbf{(iii)} We put $\Omega := \{ x \in [a_1,1] \, : \, \vphi (x,z;B) \in e^{\i \xi} \R_+ \}$ and assume that
$\Omega$ is infinite. Then, as before in Lemma \ref{l phi values},  $\Omega$ has a limit point $x_0 \in \Omega$,
and combining this with statement (ii) of the lemma, one can show that $x_0 = a_1$.

Since $B \equiv 0$ on $[0,a_1]$,  $\vphi$ satisfies
$ \vphi (x, z;B)   =  1 - z^2 \int_{a_1}^x (x-s) B (s) y (s) \d s $ (in particular, $\vphi (a_1, z;B) =1$).
This and (\ref{e Bneq0 a1x0}) implies that $x_0 \neq a_1$
in the same way as before.
\end{proof}

\begin{proof}[Proof of Lemma \ref{l 0 int} in the case $b_1 = 0$]
The proof requires only the change of Lemma \ref{l phi values} to Lemma \ref{l phi values b1=0}.
Indeed, one can see that $\Omega := \{ x \in (0,1) \, : \, b_1+\ep_1 < B_0 (x) < b_2 -\ep_1\}$ is a subset of
$(a_1,1]$. So Lemma \ref{l phi values b1=0} can be applied to show the existence of the sets $\Omega_1$ and
$\Omega_2$.
\end{proof}

Now \textbf{the proof of Corollary \ref{c k nonmin} is complete}.

The following changes are needed in connection with Proposition \ref{p semi-c}.
As before, for a suitable branch of the argument function,
$ \arg \vphi^2 (x,z;B_0) $ is differentiable.
However, $\pa_x \arg \vphi^2 (x,z;B_0) \neq 0$ only for $x\in (a_1,1]$
 (this follows from Lemma \ref{l phi values b1=0} (ii)).

\begin{prop} \label{p semi-c b1=0}
Let $B_0 \in \Ext $ have a switch point $x_0 \in (a_1,1)$.
Then all the assertions of Proposition \ref{p semi-c} hold true.
\end{prop}

Since $x_0 \in (a_1,1)$, the proof of Proposition \ref{p semi-c} works without changes.
Note that, for $B_0 \in \Ext $, the definition of $a_1$ implies that $B_0$ has \emph{no} switch points
in $[0,a_1)$.

After all these modifications, \textbf{the proofs of Proposition \ref{p singsw int0} and Corollary \ref{c
B in Exts}} requires no changes for the case $b_1 =0$ except the use of Proposition \ref{p semi-c b1=0}. Now
\textbf{the proof of Theorem \ref{t B in Exts} for $\alpha \neq 0$
is complete.}

\subsection{The case $\alpha = 0$. \label{ss case a=0}}

In this subsection we will study quasi-eigenvalues on the ray $\i \R_+$.
In this case, the analysis is simpler since
\begin{equation} \label{e vphipsi real}
\vphi (x,z;B) \text{ and } \psi (x,z;B) \text{ are real when } \ \ z \in \i \R_+
\text{ and } B \in L^1_\R (0,1).
\end{equation}

\begin{lem} \label{l Delta pm}
Let $B_0 \in \Ad $ and $\k_0 \in \i \R_+ \cap K (B_0)$.
If $B_0 \not \equiv b_1$ and $B_0 \not \equiv b_2$, then there exist
$B_{\Delta,+}, B_{\Delta,-} \in \Ad - B_0$ such that $[\pa_{B} F (\k_0,B_0)] (B_{\Delta,+}) >0$ and $[\pa_{B} F (\k_0,B_0)] (B_{\Delta,-}) <0$.
\end{lem}

\begin{proof}
It follows from (\ref{e phi psi alt}), (\ref{e vphipsi real}), and $\k_0 = \i \beta$, $\beta>0$, that
\begin{equation} \label{e kpsi neq 0}
\k_0 \left[- \k_0 \psi (1,\k_0;B_0) + \i \pa_x \psi (1,\k_0;B_0) \right] \in \R\setminus\{ 0 \} .
\end{equation}
The assumptions of the lemma imply that there exist sets $\Omega_1, \Omega_2 \subset [0,1]$ of positive measure and $\epsilon>0$ such that $B_0 (x) > b_1+\epsilon$ for $x \in \Omega_1$ and $B_0 (x) < b_2-\epsilon$ for $x \in \Omega_2$.
Put $B_{\Delta,j} = (-1)^j \epsilon \chi_{\Omega_j}$, $j=1,2$.
Then, $B_{\Delta,j} \in \Ad - B_0$, $j=1,2$.
Since the real continuous function $\vphi (\cdot,\k_0;B_0)$ has at most finite number of zeroes on $[0,1]$, we see that
\[
\int_0^1
\vphi^2 (s,\k_0;B_0) \ B_{\Delta,1} (s) \ \d s < 0 , \ \
\int_0^1
\vphi^2 (s,\k_0;B_0) \ B_{\Delta,2} (s) \ \d s > 0 .
\]
Combining this with (\ref{e paBD F}) and (\ref{e kpsi neq 0}), one gets the statement of the lemma.
\end{proof}

\begin{prop} \label{p k nonmin alpha0}
Let $B_0 \in \Ad $ and $\k_0 \in \i \R_+ \cap K (B_0)$.
If $B_0 \not \equiv b_1$ and $B_0 \not \equiv b_2$, then
there exists $\beta_1>0$ such that $\k_0 - \i \beta_1  \in K (\Ad)$.
\end{prop}

\begin{proof}
We give a proof for the case when $\k_0$ is a quasi-eigenvalue of multiplicity $r \geq 2$.
The arguments for the case of a simple quasi-eigenvalue are simpler in details.

So $\pa_z^r F (\k_0; B_0) \neq 0$.
It follows from (\ref{e vphipsi real}) that $F (\i \beta; B_0) \in \R$ when $\beta \in \R_+$.
Hence, $\pa_z^r F (\k_0; B_0) \in \i^r \R \setminus \{0\}$.
By Lemma \ref{l Delta pm}, we can choose $B_\Delta \in \Ad - B_0$ such that
\[
\arg \left( -\frac{r! \, [\pa_B F ( \k_0 ; B_0 )] (B_\Delta)}{\pa_z^r F (\k_0; B_0)} \right) = -r \pi/2 \ ( \modn 2 \pi ) .
\]
So one can choose $\sqrt[\r]{\cdot}$ such that
\begin{equation} \label{e arg sqrt h}
\arg \sqrt[\r]{ -\frac{r! \, [\pa_B F ( \k_0 ; B_0 )] (B_\Delta)}{\pa_z^r F (\k_0; B_0)} } = - \pi/2 \ (\modn 2 \pi ).
\end{equation}
Applying Proposition \ref{p per k mult}, we see that for $\zeta>0$
small enough, the $r$ branches of the Puiseux series (\ref{e k
Pser}) give all the quasi-eigenvalues of $B_0 + \zeta B_\Delta$ that
tend to $\k_0$ as $\zeta \to 0$. Let us choose in (\ref{e k Pser})
the branch of $(\cdot)^{1/r}$ analytic on $\R_+$ and such that
$1^{1/r} = 1$. Then it follows from (\ref{e arg sqrt h}) that there
exist a quasi-eigenvalue $k_1 (\zeta)$ of $B_0 + \zeta B_\Delta$,
$\zeta > 0$, with asymptotics $k_1 (\zeta) = \k_0 - \i c
|\zeta|^{1/r} + o (|\zeta|^{1/r}) $, $c >0$. The other branches of
$k (\zeta)$ have the same asymptotics for $\zeta >0$ with constants
$c \not \in \R_+$. Since $K (B_0 + \zeta B_\Delta)$ is symmetric
w.r.t. $\i \R$, $k_1 (\zeta) $ stays on
$\i \R_+$ for $\zeta > 0$ small enough.

Summarizing, we see that the structures $B_0 + \zeta B_\Delta$ belong to $\Ad$ for $\zeta \in [0,1]$,
and that for $\zeta > 0$ small enough, one of quasi-eigenvalues of $B_0 + \zeta B_\Delta$ may be written in the form
$\k_0 - \i \beta (\zeta)$ with $\beta (\zeta) > 0$.
\end{proof}

\begin{proof}[Proof of Theorem \ref{t J alpha0}.]
Let $\k_0 = \i \I (0) \in \i \R_+$ be a quasi-eigenvalue of $B_0 \in \Ad$.
Then Proposition \ref{p k nonmin alpha0} implies that either $B_0 \equiv b_1$
or $B_0 \equiv b_2$. Proposition \ref{p const struct} completes the proof.
\end{proof}

\section{The proof of Theorem \ref{t nonlin}: nonlinear eigenvalues and restrictions on switch points.}
\label{s sw pt}


Let $B$ be an optimal structure for a frequency $\alpha \neq 0 $ and
let $\k =\alpha + \i \I (\alpha) \in K (B)$ be a corresponding
optimal quasi-eigenvalue. By Theorem \ref{t B in Exts}, $B$ is a
piecewise constant function taking only the values $b_1$ and $b_2$.
In other words, $B$ has at most finite number of switch points $\{
x_j \}_{j=1}^n$ where $B$ changes its value either from $b_1$ to
$b_2$ or inversely from $b_2$ to $b_1$. Note that the values of $B$
at the switch points and the endpoints $x=0$ and $x=1$ are not
important for the quasi-eigenvalue problem.


We will use the notation of Proposition \ref{p semi-c} with
the continuous in $x$ branch $\xi (x)$ of the multifunction
$\arg \vphi^2 (x,\k;B)$ fixed by $\xi (0) = 0$.

\begin{lem} \label{l var xi}
If an optimal structure $B$ is constant on an interval $(\wt x_1 , \wt x_2 ) $, then
$|\xi (\wt x_1) - \xi (\wt x_2)| \leq \pi$.
\end{lem}

\begin{proof}
Assume the contrary, i.e., $ | \xi (\wt x_1) - \xi (\wt x_2) | > \pi $.
Then (\ref{e paBD F}) implies that the set of quasi-arguments $A
(\k;B)$ defined by (\ref{e Aset}) contains an arc of length
greater than $\pi$. So $0$ is an interior point the convex set $
[\pa_B F (\k;B)] (\Ad - B) $. Proposition \ref{p k nonmin}
implies that $\k$ is not an optimal quasi-eigenvalue, a
contradiction.
\end{proof}

Consider first the non-degenerate case $0<b_1 < b_2$.
Since $\k \in \C_+ \setminus \i \R$, the function $\xi $ is continuously differentiable on $[0,1]$ and
$\xi' (x) \neq 0$ for $x \in (0,1]$. So $\xi' (x) $ keeps its sign on
$(0,1]$.
In fact, for $x \in (0,1]$,
\begin{equation} \label{e xi pr}
\xi' (x) < 0 \ \text{ if } \ \re \k > 0, \ \  \text{ and } \ \ \xi' (x) > 0 \ \text{ if } \ \re \k < 0 .
\end{equation}
(see (\ref{e arg vphi2}) and the proof of statement (iii) of Lemma \ref{l phi values} for details).

\begin{thm} \label{t sw point}
Let $0<b_1 < b_2$. Let $B$ be an optimal structure for a frequency
$\alpha \in \R\setminus\{0\}$ and let $\k $ be a corresponding
optimal quasi-eigenvalue.  Then there exists $\omega \in [-\pi,\pi)$
such that, on the interval $(0,1)$, $B (x)$ changes its value from
$b_1$ to $b_2$ exactly when $\vphi^2 (x,\k ;B)$ intersects the ray
$e^{\i \omega} \R_+$ and $B$ changes its value from $b_2$ to $b_1$
exactly when $\vphi^2 (x,\k ; B)$ intersects the ray $e^{\i \omega }
\R_-$.
\end{thm}

\begin{proof}
We assume that $\xi'(x) < 0$ on $(0,1]$ (arguments for the case $\xi'(x) > 0$ are similar).

\emph{Case 1.} Assume, first, that $B$ is constant on $(0,1)$. Then Lemma \ref{l var xi} implies that
$\xi (1) \geq -\pi $ and the statement of the theorem is valid both with
$\omega = -\pi $ and $\omega = 0 $.

\emph{Case 2.} Assume that $B$ has only one switch point $x_1$.
Then, in the case when $B$ changes its value from $b_1$ to $b_2$
at $x_1$,  Lemma \ref{l var xi} yields the statement of the theorem with
$\omega = \xi (x_1) $, and in the opposite case with $\omega = \xi (x_1) + \pi $.

\emph{Case 3.} Consider the case when there are at least two switch points.
Let the set  $\{ x_j \}_{j=1}^n$ of switch points be
naturally ordered $ 0 < x_1 < x_2 < \dots < x_{n-1} < x_n < 1 $. To
be specific, assume that $B$ changes its value from $b_1$ to $b_2$
at the first switch point $x_1$.

By Lemma \ref{l var xi}, $\xi (x_1) \in [-\pi, 0)$.
Assign $\omega = \xi (x_1) $.
From Proposition \ref{p semi-c} and the assumption $\xi'(x_1) < 0$, we see that
the set $A (\k;B) $ contains the semi-circles
 $ \{ e^{\i[\xi (x_j) +\xi_1+s]} \, : \, s \in (0,\pi) \} $ for odd $j$
 and the semi-circles $ \{ e^{\i[\xi (x_j) +\xi_1+s]} \, : \, s \in (-\pi,0) \} $
 for even $j$. Since $\k$ is an optimal
quasi-eigenvalue, all these semi-circles coincide (see the proofs of
Propositions \ref{p singsw int0} and Corollary \ref{c B in Exts}).
This means that
\begin{eqnarray}
\xi (x_j) & = & \xi (x_1) \ \ (\modn 2 \pi)\ \  \text{ for odd } \ \ j , \\
\xi (x_j) & = & \xi (x_1) + \pi \ \ (\modn 2 \pi) \  \text{ for even } \ \ j .
\end{eqnarray}
Thus, $B$ may change its value
from $b_1$ to $b_2$ only at the points where $\xi (x) = \omega \ ( \modn 2 \pi)$ and
from $b_2$ to $b_1$ only where $\xi (x) =  \omega + \pi \ ( \modn 2 \pi)$.
On the other side, Lemma \ref{l var xi} implies that $B$ indeed changes its value
each time when $\xi (x) =  \omega \ (\modn 2 \pi)$ or $\xi (x) =  \omega + \pi \ (\modn 2 \pi)$.
This completes the proof.

In the case when $B$ changes its value from $b_2$ to $b_1$
at $x_1$, analogous arguments produce $\omega = \xi (x_1) + \pi$.
\end{proof}

\begin{proof}[The proof of Theorem \ref{t nonlin} in the case $b_1 > 0$.]
Let $\alpha=0$. Then an optimal quasi-eigenvalue $\k$ exists only if $b_2 >1$.
The optimal structure is $B\equiv b_2$, and
the corresponding mode $\vphi ( \cdot,\k;b_2)$ is a real-valued function with a finite number of zeroes (see Section \ref{ss case a=0}).
Hence, for arbitrary $\theta \in (0,\pi/2)$, the function $y = e^{\i \theta} \vphi$
is a solution of the nonlinear eigenvalue problem (\ref{e nonlin eq}), (\ref{e bc0}), (\ref{e bc1}).

Now consider $\alpha \neq 0$ and take $\omega$ as in Theorem
\ref{t sw point}. Then, to obtain a solution of the nonlinear problem
(\ref{e nonlin eq}), (\ref{e bc0}), (\ref{e bc1}), one can put $y
(\cdot) = e^{\i (\pi -\omega)/2} \vphi ( \cdot,\k;B) $ in the case
$\re \k > 0$, and $y (\cdot) = e^{- \i \omega/2} \vphi ( \cdot,\k;B)
$ in the case $\re \k < 0$. This fact follows immediately from
(\ref{e xi pr}) and Theorem \ref{t sw point}.
\end{proof}

If $b_1 = 0$, some technical complications arise
since (\ref{e xi pr}) is valid only on $(a_1,1]$, while $\xi' (x) = 0$ for $x \in [0,a_1]$.
Recall that $[0,a_1]$ is the greatest interval of the form $[0,x_0]$ such that
$B (x) =0$ a.e. on $[0,x_0]$, see Section \ref{ss case b1=0} for details.

\begin{proof}[The proof of Theorem \ref{t nonlin} in the degenerate case $b_1 = 0$.]
In the cases when $\alpha = 0$ or $a_1 =0$,
the proof is the same as in the case $b_1>0$.

Consider the remaining possibility when $\alpha \neq 0$ and $B (x) = 0$ on $(0,a_1)$ with $0<a_1<1$.
Then $a_1 $ is the first switch point of $B$ and $\vphi (x,\k;B) = 1$ for $x \in [0,a_1]$.
With no loss of generality, assume additionally that $\xi'(x) < 0$ on $(a_1,1]$
(that is, we assume $\re \k > 0$).

Let us show that, on the interval $(a_1, 1)$, the optimal structure
$B$ changes its value from $0$ to $b_2$ exactly when $\vphi^2 (x,\k_0 ;B)$ intersects
$\R_+$ and from $b_2$ to $0$ exactly when $\vphi^2 (x,\k_0 ; B)$ intersects $\R_-$.


Indeed, under the above assumptions, the arguments of Proposition \ref{p semi-c}
applied to the first switch point $x_1 = a_1$
imply that $ A(\k; B) $ contains $\{ e^{\i[\xi_1+s]} \, : \, s \in (0,\pi)\}$ (note that $\xi (a_1) = 0$).
Since (\ref{e xi pr}) holds true on $(a_1,1]$,  we
can deal with all the other switch points $x_2$, $x_3$, $\dots$ in the same way as
in the proof of Theorem \ref{t sw point}. As a result, we obtain
that the semi-circle $\{ e^{i[\xi_1+s]} \, : \, s \in (0,\pi)\}$ generated by the first switch point and
all the semi-circles
\begin{eqnarray}
 \{ e^{\i[\xi (x_j) +\xi_1+s]} \, : \, s \in (0,\pi) \} \ \ \text{ for odd }  \ j > 1, \ \ \
 \{ e^{\i[\xi (x_j) +\xi_1+s]} \, : \, s \in (-\pi,0) \} \ \ \text{ for even } \ j ,
\end{eqnarray}
generated by the other switch points coincide.
This yields the desired statement.

Now one can take $y (x) = e^{-\i \pi /2 } \vphi (x , \k ; B) $ and check that $y$ is a solution of
the nonlinear problem (\ref{e nonlin eq}), (\ref{e bc0}), (\ref{e bc1}).

Note that in the case $\xi'(x) > 0$ on $(a_1,1]$
(i.e., when $\re \k <  0$), these arguments produce $y (x) =  \vphi (x , \k ; B) $.
\end{proof}

\emph{Acknowledgements.} The author is grateful to Richard Froese for bringing this problem to his attention and for a number of interesting discussions.
This work was financially supported by Elena Braverman, Richard Froese, and
the Pacific Institute for the Mathematical Sciences.
The author would like to thank Elena Braverman for the hospitality of the University of Calgary, and Richard Froese for the hospitality of the University of British Columbia.


\begin{thebibliography}{38}

\bibitem{AAD01} A.A. Abramov, A. Aslanyan, E.B. Davies, Bounds on complex eigenvalues and resonances,
J. Phys. A 34 (2001), no.1, 57--72.

\bibitem{A98} S. Agmon, A perturbation theory of resonances, Comm. Pure Appl. Math. 51 (1998), no.11-12, 1255--1309.

\bibitem{AASN03} Y. Akahane, T. Asano, B. Song, S. Noda, High-Q photonic
nanocavity in a two-dimensional photonic crystal, Nature 425 (2003)
944--947.

\bibitem{A75} D.Z. Arov, The realization of a canonical system with dissipative boundary conditions at
one end of a segment in terms of the coefficient of dynamic flexibility,
Sibirsk. Mat. Zh. 16 (1975), no.3, 440--463; English transl.: Siberian Math. J. 16 (1975), no.3, 335--352.


\bibitem{BOYa04} M. Burger, S. Osher, E. Yablonovitch, Inverse problem
techniques for the design of photonic crystals. IEICE Trans.
Electron. 87 (2004) 258--265.

\bibitem{CMcL90} S.J. Cox, J.R. McLaughlin,
Extremal eigenvalue problems for composite membranes. I, II,
Appl. Math. Optim. 22 (1990), no. 2, 153--167, 169--187.


\bibitem{CZ95}
S. Cox, E. Zuazua, The rate at which energy decays in a string damped at one end, Indiana Univ. Math. J. 44 (1995), no.2,
545--573.

\bibitem{DM76} H. Dym, H.P. McKean, Gaussian Processes, Function Theory, and the Inverse Spectral Problem,
Academic Press, New York, 1976.

\bibitem{F97} R. Froese, Asymptotic distribution of resonances in one dimension, J. Differential Equations 137 (1997),
no.2, 251--272.

\bibitem{GK71} T.M. Gataullin, M.V. Karasev,
On the perturbation of the quasilevels of a Schr\"odinger operator with complex potential, Teoret. Mat. Fiz. 9 (1971),
no.2, 252--263 (Russian); English transl.: Theoretical and Mathematical Physics
9 (1971), no.2, 1117--1126.

\bibitem{GP97} G.M. Gubreev, V.N. Pivovarchik,
Spectral analysis of the Regge problem with parameters,  Funktsional. Anal. i Prilozhen. 31 (1997), no.1, 70--74 (Russian);
English transl.: Funct. Anal. Appl. 31 (1997), no. 1, 54--57.


\bibitem{HBKW08} P. Heider, D. Berebichez, R.V. Kohn, M.I. Weinstein,
Optimization of scattering resonances, Struct. Multidisc. Optim. 36 (2008) 443--456.

\bibitem{JJWM08} J.D. Joannopoulos, S.G. Johnson, J.N. Winn, R.D. Meade, Photonic
Crystals: Molding the Flow of Light, Princeton University Press,
2008.

\bibitem{KK68_II}
I.S. Kac and M.G. Krein, On the spectral functions of the string,
Supplement II in Atkinson, F. Discrete and continuous boundary
problems. Mir, Moscow 1968. Engl. transl.:
 Amer. Math. Soc. Transl., Ser. 2, 103 (1974) 19--102.


\bibitem{KS08}
 C.-Y. Kao, F. Santosa, Maximization of the quality factor of an optical resonator, Wave Motion 45 (2008) 412--427.


\bibitem{Ka12_optKN} I.M. Karabash, Optimization of quasi-normal eigenvalues for
Krein-Nudelman strings, in preparation.


\bibitem{K51}
M.G. Krein, On certain problems on the maximum and minimum of characteristic values and on the Lyapunov zones of stability,
Prikl. Mat. Meh. 15 (1951) 323--348 (Russian); English transl.: Amer. Math. Soc. Transl. (2) 1 (1955) 163--187.

\bibitem{K93} M.G. Krein, Selected works. I. Complex analysis, extrapolation, interpolation,
Hermitian-positive functions and related topics, Akad. Nauk Ukrainy, Kiev, 1993 (Russian).

\bibitem{KN79}
M.G. Krein, A.A. Nudelman, On direct and inverse problems for the boundary
dissipation frequencies of a nonuniform string, Dokl. Akad. Nauk SSSR 247 (1979), no. 5, 1046--1049 (Russian);
English transl.: Soviet Math. Dokl. 20 (1979), no.4, 838--841.

\bibitem{KN81} M.G. Krein, A.A. Nudelman,
Representation of entire functions that are positive on the real axis, or on the half axis, or outside a finite interval,
Mat. Issled. 61 (1981), 40--59 (Russian).

\bibitem{KN89}
M.G. Krein, A.A. Nudelman, Some spectral properties of a nonhomogeneous
string with a dissipative boundary condition, J. Operator Theory 22 (1989) 369--395 (Russian).

\bibitem{LSV03} R.P. Lipton, S.P. Shipman, S. Venakides, Optimization of resonances
of photonic crystal slabs,  in: Proceedings SPIE, vol. 5184, pp. 168--177.

\bibitem{N69} M.A. Naimark, Linear differential operators, second ed., Nauka, Moscow, 1969 (Russian);
English transl.: Parts I, II. Frederick Ungar Publishing Co., New York, 1967-68.

\bibitem{PT87} J. P\"{o}schel, E. Trubowitz,
Inverse spectral theory,  Pure and Applied Mathematics 130, Academic Press, Boston, 1987.

\bibitem{RS78IV} M. Reed, B. Simon, Methods of modern mathematical physics. IV. Analysis of operators,
 Academic Press , New York-London, 1978.

\bibitem{Sh07} A.A. Shkalikov, Spectral analysis of the Redge problem. J. Math.
Sci. 144 (2007), no. 4, 4292--4300.


\end{thebibliography}
\end{document}